\documentclass[11pt,reqno]{amsart}
\pdfoutput=1 

\usepackage[english]{babel}
\usepackage[utf8]{inputenc}
\usepackage[T1]{fontenc}
\usepackage{lmodern} \normalfont 
\DeclareFontShape{T1}{lmr}{bx}{sc} { <-> ssub * cmr/bx/sc }{}
\usepackage{microtype}	

\usepackage{amssymb}
\usepackage{romannum}
\usepackage{amsmath}
\usepackage{comment}
\usepackage{amsthm}
\usepackage{accents}
\usepackage{algorithm}
\usepackage{algpseudocode}
\usepackage{mathtools}
\usepackage{mathrsfs}
\usepackage{accents} 
\usepackage{subdepth} 
\usepackage{etoolbox}
\usepackage{siunitx}
\usepackage{dsfont} 
\usepackage{bm}

\usepackage{tabularx}

\usepackage{graphicx}
\usepackage{color}
\usepackage[table,xcdraw,dvipsnames]{xcolor}
\usepackage{circuitikz}
\usepackage{tikz}
\usetikzlibrary{calc,positioning,shapes}
\usetikzlibrary{patterns,decorations.pathmorphing,decorations.markings}
\usetikzlibrary{external}
\usepackage{pgfplots}
\pgfplotsset{compat=newest}
\usepackage[margin=10pt,font=small,labelfont=bf,labelsep=endash]{caption}
\usepackage{subcaption}

\usepackage{booktabs}
\usepackage{paralist}
\usepackage{soul}
\usepackage{listings}

\numberwithin{equation}{section}

\usepackage{algorithm}
\usepackage{algorithmicx}
\usepackage{algpseudocode}

\usepackage{cite}
\usepackage{multirow} 
\usepackage{enumitem} 
\setlist[enumerate]{label=(\roman*)}
\usepackage{xspace}

\usepackage[colorlinks,linkcolor=teal,citecolor=teal,urlcolor=teal]{hyperref}
\usepackage[nameinlink,noabbrev]{cleveref}

\usepackage{booktabs} 

\theoremstyle{plain}
\newtheorem{theorem}{Theorem}[section]
\newtheorem{proposition}[theorem]{Proposition}
\AddToHook{env/proposition/begin}{\crefalias{theorem}{proposition}}
\Crefname{proposition}{Proposition}{Propositions}

\newtheorem{lemma}[theorem]{Lemma}
\AddToHook{env/lemma/begin}{\crefalias{theorem}{lemma}}
\Crefname{lemma}{Lemma}{Lemmata}

\AddToHook{env/corollary/begin}{\crefalias{theorem}{corollary}}
\Crefname{corollary}{Corollary}{Corollaries}

\newtheorem{remark}[theorem]{Remark}
\AddToHook{env/remark/begin}{\crefalias{theorem}{remark}}
\Crefname{remark}{Remark}{Remarks}

\newtheorem{definition}[theorem]{Definition}
\AddToHook{env/definition/begin}{\crefalias{theorem}{definition}}
\Crefname{definition}{Definition}{Definitions}

\newtheorem{assumption}[theorem]{Assumption}
\AddToHook{env/assumption/begin}{\crefalias{theorem}{assumption}}
\Crefname{assumption}{Assumption}{Assumptions}

\AddToHook{env/example/begin}{\crefalias{theorem}{example}}
\Crefname{example}{Example}{Examples}

\AddToHook{env/problem/begin}{\crefalias{theorem}{problem}}
\Crefname{problem}{Problem}{Problems}

\usepackage{ifthen}
\newboolean{preprint}
\setboolean{preprint}{false} 	

\newcommand{\N}{\ensuremath\mathbb{N}}

\newcommand{\R}{\ensuremath\mathbb{R}}

\newcommand{\smoothFunctions}[3][]{\ifthenelse{\equal{#1}{}}{\mathcal{C}^{#2}}{\mathcal{C}_{#1}^{#2}}(#3)}

\newcommand{\dist}[1][]{\ifthenelse{\equal{#1}{}}{\mathbb{D}}{#1_{\mathbb{D}}}}

\newcommand{\Spsd}[1]{\mathcal{S}^{#1}_{\succeq}}

\newcommand{\Frob}{\mathrm{F}}

\newcommand{\T}{\ensuremath\mathsf{T}}

\newcommand{\integrate}[1]{\mathrm{d}#1}

\newcommand{\ds}{\integrate{s}}
\newcommand{\dt}{\integrate{t}}
\newcommand{\dz}{\integrate{z}}
\newcommand{\dtheta}{\integrate{\theta}}

\DeclareMathOperator{\ee}{e}

\DeclareMathOperator{\real}{Re}

\DeclareMathOperator{\spann}{span}

\DeclareMathOperator{\vec2}{vec}

\newcommand{\calH}{\mathcal{H}}

\newcommand{\calJ}{\mathcal{J}}
\newcommand{\calK}{\mathcal{K}}
\newcommand{\calL}{\mathcal{L}}

\newcommand{\calO}{\mathcal{O}}
\newcommand{\calP}{\mathcal{P}}
\newcommand{\calQ}{\mathcal{Q}}
\newcommand{\calR}{\mathcal{R}}
\newcommand{\calS}{\mathcal{S}}
\newcommand{\calT}{\mathcal{T}}
\newcommand{\calU}{\mathcal{U}}
\newcommand{\calV}{\mathcal{V}}

\newcommand{\stx}{\ensuremath{\bm{x}}}
\newcommand{\state}{\stx}
\newcommand{\stateDim}{n}
\newcommand{\stateRed}{\reduce{\state}}
\newcommand{\stateDimRed}{r}
\newcommand{\stateBal}{\bar{\state}}
\newcommand{\inp}{\ensuremath{\bm{u}}}
\newcommand{\inpDim}{m}
\newcommand{\out}{\ensuremath{\bm{y}}}
\newcommand{\outDim}{p}

\newcommand{\fw}{\ensuremath{\bm{w}}}

\newcommand{\fz}{\ensuremath{\bm{z}}}

\newcommand{\fA}{\ensuremath{\bm{A}}}
\newcommand{\fB}{\ensuremath{\bm{B}}}
\newcommand{\fC}{\ensuremath{\bm{C}}}
\newcommand{\fD}{\ensuremath{\bm{D}}}
\newcommand{\fE}{\ensuremath{\bm{E}}}

\newcommand{\fG}{\ensuremath{\bm{G}}}

\newcommand{\fI}{\ensuremath{\bm{I}}}
\newcommand{\fJ}{\ensuremath{\bm{J}}}

\newcommand{\fM}{\ensuremath{\bm{M}}}
\newcommand{\fN}{\ensuremath{\bm{N}}}

\newcommand{\fP}{\ensuremath{\bm{P}}}
\newcommand{\fQ}{\ensuremath{\bm{Q}}}
\newcommand{\fR}{\ensuremath{\bm{R}}}
\newcommand{\fS}{\ensuremath{\bm{S}}}

\newcommand{\fU}{\ensuremath{\bm{U}}}
\newcommand{\fV}{\ensuremath{\bm{V}}}
\newcommand{\fW}{\ensuremath{\bm{W}}}
\newcommand{\fX}{\ensuremath{\bm{X}}}
\newcommand{\fY}{\ensuremath{\bm{Y}}}
\newcommand{\fZ}{\ensuremath{\bm{Z}}}

\newcommand{\fmu}{\ensuremath{\bm{\mu}}}

\newcommand{\fDelta}{\ensuremath{\bm{\Delta}}}

\newcommand{\fSigma}{\ensuremath{\bm{\Sigma}}}

\newcommand{\IdMat}{\mathbf{I}}

\newcommand{\fcP}{\mathbf{\calP}}
\newcommand{\fcQ}{\mathbf{\calQ}}
\newcommand{\zeroVec}{\mathbf{0}}
\newcommand{\zeroMat}{\mathbf{0}}

\newcommand{\bfX}{\mathbf{x}}

\newcommand{\switch}{q}
\newcommand{\switchingSet}{\calJ}

\newcommand{\tol}{\texttt{tol}}
\newcommand{\ir}{\alpha}

\newcommand{\system}{\Sigma}
\newcommand{\reduce}[1]{\tilde{#1}}
\newcommand{\reduceBis}[1]{\hat{#1}}
\newcommand{\systemRed}{\reduce{\system}}
\newcommand{\switchedSys}{\system_{\switch}}
\newcommand{\switchedSysRed}{\systemRed_{\switch}}
\newcommand{\switchedSysRedJumps}{\hat \system_{\switch}}
\newcommand{\switchedSysRedJumpsPert}{\hat \system_{\switch,\Delta}}

\newcommand{\gleMat}[1]{\mathbf{\mathscr{#1}}}

\newcommand{\gleG}{\gleMat{G}}
\newcommand{\gleH}{\gleMat{H}}

\newcommand{\gleM}{\gleMat{M}}
\newcommand{\gleN}{\gleMat{N}}
\newcommand{\lyapOper}{\calL}
\newcommand{\gleOper}{\Pi}
\newcommand{\perturbation}{\mu}

\newcommand{\indKS}{\ell}

\newcommand{\lineWidth}{1.2pt}
\newcommand{\imageWidth}{2.0in}
\newcommand{\imageHeight}{1.8in}
\definecolor{mycolor1}{rgb}{0.00000,0.44700,0.74100}
\definecolor{mycolor2}{rgb}{0.85000,0.32500,0.09800}
\definecolor{mycolor3}{rgb}{0.92900,0.69400,0.12500}
\definecolor{mycolor4}{rgb}{0.46600,0.67400,0.18800}
\definecolor{mycolor5}{rgb}{0.49400,0.18400,0.55600}
\definecolor{mycolor6}{rgb}{0.30100,0.74500,0.93300}

\newcommand{\abbr}[1]{\textsf{#1}\xspace}
\newcommand{\FOM}{\abbr{FOM}}
\newcommand{\GLE}{\abbr{GLE}}
\newcommand{\LE}{\abbr{LE}}
\newcommand{\GLEs}{\abbr{GLEs}}
\newcommand{\LMI}{\abbr{LMI}}
\newcommand{\LMIs}{\abbr{LMIs}}
\newcommand{\ROM}{\abbr{ROM}}
\newcommand{\ROMs}{\abbr{ROMs}}
\newcommand{\MOR}{\abbr{MOR}}

\newcommand{\ODE}{\abbr{ODE}}
\newcommand{\SLS}{\abbr{SLS}}
\newcommand{\ODEs}{\abbr{ODEs}}
\newcommand{\PDE}{\abbr{PDE}}
\newcommand{\PDEs}{\abbr{PDEs}}

\newcommand{\PBR}{\abbr{PBR}}

\newcommand{\switchedSysRedJumpsGLE}{\hat \system_{\switch,\GLE}}
\newcommand{\switchedSysRedJumpsLE}{\hat \system_{\switch,\LE}}
\title[Model Reduction for SLS via Generalized Lyapunov Equations]{Model Reduction for Switched Linear Systems via Generalized Lyapunov Equations}
\author[M.~Manucci \and B.~Unger]{Mattia Manucci${}^\star$ \and Benjamin Unger${}^\star$}

\address{${}^{\star}$ Institute for Applied and Numerical Mathematics, Karlsruhe Institute of Technology, 76131 Karlsruhe, Germany}
\email{\{mattia.manucci,benjamin.unger\}@kit.edu}

\newcommand{\ourKeywords}{generalized Lyapunov equation, model order reduction, switched systems, balanced truncation, error bound}
\date{\today}

\begin{document}
	
\begin{abstract}
In this work, we study projection-based model order reduction (MOR) for switched linear systems (SLS) in control form, where the projection matrices are obtained from the solutions of generalized Lyapunov equations (GLEs). We investigate how numerical inaccuracies in solving the GLEs propagate through the MOR process and impact the accuracy and reliability of the resulting reduced-order model. This highlights the importance of accounting for such inaccuracies, motivating the introduction of a novel error bound to quantify and control the error in the approximation of the GLE solution.

Moreover, classical balanced truncation error estimates for SLS are neither theoretically sound nor practically applicable, as they rely on restrictive assumptions requiring several linear matrix inequalities (LMIs) to be satisfied exactly by numerically computed GLE solutions. To address these limitations, we propose a new MOR framework for SLS, termed piecewise balanced reduction (PBR). The approach is based on solving multiple GLEs and constructing projection matrices that are piecewise constant in time. By extending the standard balanced truncation error bound for SLS, we show that the PBR framework effectively controls errors arising from inexact LMI satisfaction. In addition, the proposed error bound captures the influence of the piecewise constant in time projection matrices. Altogether, this makes the PBR approach applicable to a broad and flexible class of switched linear systems. Numerical experiments are presented to support the theoretical results.
\end{abstract}
\maketitle
	
%
%
	
\section{Introduction}
We consider \emph{switched linear system} (\SLS) of the form 
\begin{equation}
	\label{eqn:sDAE}
	\switchedSys \quad \left\{\quad \begin{aligned}
	 \dot{\stx}(t)  &= \fA_{\switch(t)} \stx(t) +\fB_{\switch(t)}\inp(t), & \stx(t_0) &= \zeroVec, \\
		\out(t) &= \fC_{\switch(t)}\stx(t),\\
	\end{aligned}\right.
\end{equation}
where the symbols $\stx(t)\in\R^{\stateDim}$, $\inp(t)\in\R^{\inpDim}$, and $\out(t)\in\R^{\outDim}$ denote the \emph{state}, the controlled \emph{input}, and the measured \emph{output}, respectively. In~\eqref{eqn:sDAE}, $\switch\colon \R\to \switchingSet\vcentcolon=\{1,\ldots,M\}$ is the external switching signal, which we assume to be an element of the set of allowed switching signals
\begin{equation}
	\label{eqn:suitableSwitchingSignals}
	\calS \vcentcolon= \{\switch\colon \R\to \switchingSet \mid \switch \text{ is right continuous with a locally finite number of jumps}\}.
\end{equation}
Finally, the system matrices $\fA_{j}\in\R^{\stateDim\times \stateDim}$, $\fB_j\in \R^{\stateDim\times \inpDim}$, and~$\fC_j\in\R^{\outDim\times \stateDim}$ correspond to the \emph{ordinary differential equation}~(\ODE) that is active in mode $j\in\switchingSet$. Throughout the manuscript, we assume that $\fA_j$ is Hurwitz, i.e., its spectrum lies in the open left‑half complex plane for all $j\in\switchingSet$; moreover, we assume that $\inpDim,\outDim\ll\stateDim$. In the following, we will refer to \eqref{eqn:sDAE} as the \emph{full-order model} (\FOM). Sample applications of switched systems include robot manipulators, traffic management, automatic gear shifts, and power systems; see, for example, \cite{Che04} and the references therein. 

\begin{remark}
	Note that one can generalize the switched system~\eqref{eqn:sDAE} by introducing a nonsingular symmetric positive-definite matrix $\fE_{j}\in\R^{\stateDim\times\stateDim}$ that multiplies $\dot{\stx}(t)$. In fact, our methodology naturally extends to the case $\fE_j\neq\fI_\stateDim$ by inverting $\fE_j$ and handling \(\fE_j^{-1}\fA_j\) and \(\fE_j^{-1}\fB_j\) implicitly. 
\end{remark}

If~\eqref{eqn:sDAE} has to be repeatedly evaluated, for example, in a simulation context for different inputs or switching signals, or if matrix equalities or inequalities using the system matrices $\fA_j$, $\fB_j$, $\fC_j$ in the context of synthesis have to be solved, then a large dimension $\stateDim$ of the state makes this a computationally expensive task. In such scenarios, one can rely on \MOR and replace~\eqref{eqn:sDAE} with the \emph{reduced-order model} (\ROM)
\begin{equation}
	\label{eqn:sDAE:ROM}
	\switchedSysRed \quad \left\{\quad \begin{aligned}
	 \dot{\stateRed}(t)  &= \reduce{\fA}_{\switch(t)} \stateRed(t) + \reduce{\fB}_{\switch(t)}\inp(t), &	\stateRed(t_0) &= \zeroVec, \\
		\reduce{\out}(t) &= \reduce{\fC}_{\switch(t)}\stateRed(t),\\
	\end{aligned}\right.
\end{equation}
with $ \reduce{\fA}_{j}\in\R^{\stateDimRed\times \stateDimRed}$, $\reduce{\fB}_j\in \R^{\stateDimRed\times \inpDim}$, $\reduce{\fC}_j\in\R^{\outDim\times \stateDimRed}$, and $\stateDimRed\ll \stateDim$. In many cases (see \cite{Ant05}), the reduced system matrices are obtained via the Petrov–Galerkin projection; i.e., one constructs matrices $\fV,\fW\in\R^{\stateDim\times\stateDimRed}$ that satisfy $\fW^{\T}\fV=\IdMat_{\stateDimRed}$ and then defines
\begin{align}
	\label{eqn:sDAE:ROM:matrices}
	\reduce{\fA}_j &\vcentcolon= \fW^\T \fA_j \fV, &
	\reduce{\fB}_j &\vcentcolon= \fW^\T \fB_j, &
	\reduce{\fC}_j &\vcentcolon= \fC_j \fV.
\end{align}
The goal of \MOR is thus to derive the matrices $\fV,\fW$ in a computationally efficient and robust way, such that the error $\out-\reduce{\out}$ is small in a given norm. One approach, presented in \cite{PonGB20}, is to solve a pair of \emph{generalized Lyapunov equations} (\GLEs) of the form
\begin{subequations}
	\label{eqn:GLE}
	\begin{align}
		\fA \fcP+\fcP\fA^\T+\sum_{j=1}^{M}\left(\fN_j \fcP\fN_j^\T+\fB_j\fB_j^\T \right) &= \zeroVec,\label{eqn:GLE:reach}\\
		\fA^\T \fcQ + \fcQ\fA+\sum_{j=1}^{M}\left(\fN^\T_j \fcQ\fN_j + \fC_j^\T\fC_j\right) &= \zeroVec.\label{eqn:GLE:observ}
	\end{align}
\end{subequations}
where $\fA \vcentcolon= \fA_1$ and $\fN_j \vcentcolon= \fA_j-\fA_1$ for $j=1,\ldots,M$. From the solution of \eqref{eqn:GLE}, it is possible to obtain the projection matrices $\fV,\fW$, and thus the reduced system \eqref{eqn:sDAE:ROM:matrices}. For this reason, efficiently approximating the solution of the large-scale generalized Lyapunov equation becomes crucial for \MOR. Another often overlooked aspect is that the approximation error of the \GLE solution directly influences the output discrepancy between the \FOM and the \ROM. Therefore, taking this error into account is crucial for rigorously assessing the quality of the \ROM approximation.
\subsection{Main contributions}
To deal with the large-scale setting, we apply the stationary algorithm from \cite{ShaSS16} in combination with a subspace projection framework \cite{Sim16} to solve \GLEs; this approach is widely exploited in the literature on \GLEs solution approximation. Our first contribution is the derivation of an error estimate such that, for any prescribed user tolerance $\tol$, the computed approximation $\tilde{{\fX}}$ of $\fX$,  where $\fX\in\{\fcP,\fcQ\}$ being $\fcP,\fcQ$, the exact solution of \eqref{eqn:GLE:reach} and \eqref{eqn:GLE:observ}, respectively, satisfies $\|\fX-\tilde{\fX}\|_2\le\tol$; see \Cref{teo1}. The theorem is followed by an extensive discussion that clarifies the conditions under which the error estimate can be computed efficiently.

Second, we analyze how numerical errors in the approximation of \eqref{eqn:GLE} can degrade both the quality and stability of the reduced-order model (\ROM) in \eqref{eqn:sDAE:ROM}. For a specific class of switched systems, we demonstrate how to appropriately perturb the approximated solution 
$\tilde \fX$, leveraging the error certification provided by our algorithm to ensure that the resulting \ROM remains stable while achieving arbitrarily high accuracy; see \Cref{teo4}.

The error bound for the balancing approach \cite{PonGB20} applied to \SLS relies on the assumption that the solution $\fX$ of the \GLE satisfies certain \emph{linear matrix inequalities} (\LMI). This assumption considerably restricts the number of problems in which the \MOR strategy of \cite{PonGB20} provides error guarantees. Moreover, even if the \LMI are satisfied by $\fX$, they may not be satisfied by its numerical approximation $\tilde\fX$. To overcome this limitation, we introduce a novel \emph{piecewise balancing reduction} (\PBR) approach for \SLS. After discussing the stability conditions of the associated \ROM, we derive a novel error bound that accounts for a potentially violated \LMI and a piecewise constant in time projection basis; see \Cref{sec:PMOR} and \Cref{teo6}. 

The effectiveness of our approach is verified in \Cref{sec:examples} through a synthetic example (\Cref{subsec:MSD}) and a switched system arising from the semi-discretization of a parametric \emph{partial differential equation} (\PDE) (\Cref{subsec:BlackSholes}).

\subsection{Literature review and state-of-the-art}

Several works have addressed the approximation of the \GLE solution. The existence and uniqueness of the solutions of the \GLE are established, for instance, in \cite[Thm.~3.6.1]{Dam04}, which is essentially based on the operator splitting ideas developed in \cite{Sch65}. We also mention \cite[Lem.~4.2]{Won68}, and for the more general case of generalized Sylvester equations, a similar characterization of solvability in \cite[Thm.~2.1]{JarMPR18}. The interpretation of the solution as Gramians of bilinear and stochastic linear control systems and their relation to energy functionals is discussed in \cite{BenD11}. Conditions for a fast singular value decay of the solution matrix, essential for a low-rank approximation of the solution, have been developed in \cite{BenB13,JarMPR18}. In terms of numerical schemes for the solution of the \GLE, we mention the ADI-preconditioned Krylov subspace method \cite{Dam08}, the bilinear ADI method \cite{BenB13}, the alternative linear scheme \cite{KreS15}, which can be interpreted as finding $\calH_2$ optimal search directions \cite{BenB12,BreR19}, and the stationary iteration of \cite{ShaSS16}, which we also employ here, where the inner iteration is solved with suitable increasing accuracy to minimize computational effort.

As one of our main contributions relates to the use of the \GLE for \MOR of continuous-time \SLS, we also provide an overview of the literature on \MOR for switched \ODEs. The authors of \cite{PapP14,PapP16} propose constructing \ROMs for each mode, independent of the other modes, completely ignoring the transition from one mode to another. Moreover, if a state-dependent switching signal is allowed, any approximation of the \ODE may be arbitrarily poor. Both phenomena are detailed with examples in \cite{SchU18}. In \cite{WuZ09}, an approach is proposed based on a set of coupled linear matrix inequalities, which becomes infeasible in a large-scale context. However, matrix inequalities can be used to guarantee the quadratic stability of the \ROM and to derive an error bound \cite{PetWL13}. Instead of matrix inequalities, \cite{ShaW12,GosPAF18} propose solving a set of coupled Lyapunov equations to compute Gramians, which can then be used for a balancing-based model reduction. There is no guarantee that this approach admits a solution \cite{Lib03}, and in a large-scale setting, the computational complexity may be very demanding. If the Gramians for each mode can be diagonalized simultaneously, then classical balanced truncation methods can be adapted, as discussed in \cite{MonTC12}.
In contrast, the methods reported in \cite{SchU18,PonGB20} rely on reformulating the switched system as a non-switched system by suitably interpreting the switching signal as a control input. In \cite{SchU18}, the switched system is recast as a linear system such that standard methods can be applied, while in \cite{PonGB20}, the system is recast as a bilinear system, and a balanced truncation approach for such systems is employed. In this paper, we closely follow the second strategy. Both approaches yield a priori error guarantees, provided that certain often non-general assumptions are satisfied. If prior information on the switching sequence is available, then the method in \cite{PonGB20} can be further specialized, as reported in \cite{GosPBA20}. Although we do not pursue this approach, our method can be adapted in a similar way. Interpolation-based techniques are discussed in \cite{ScaA16} for hybrid systems and in \cite{BasPWL16} for switched systems. For a data-driven approach constructing reduced models, we refer to \cite{GosPA18}. We emphasize that model reduction is closely related to realization theory. For switched systems, the associated connections are illustrated in \cite{PetG22} and the references therein. In \cite{PeiK19}, the authors discretize the control variable to obtain a switched system of autonomous \PDEs that they approximate with \ROMs in a model predictive control framework. Finally, we mention the recent work \cite{KarMUV24} where the authors deal with a certain class of optimal control problems constrained to \SLS dynamic. They employ a model predictive control strategy to address this problem, and the use of \ROM for \SLS is crucial for the efficiency of the strategy and thus its applicability in real-time scenarios.

\subsection{Organization of the manuscript}
After this introduction, we detail the \MOR method from \cite{PonGB20} in \Cref{sec:prliminaries}. In \Cref{subsec:app:GLE}, we show how to enforce error guarantees for the stationary iteration algorithm that approximates the solutions in \GLEs. Then, in \Cref{subsec:red:err}, for a certain class of \SLS, we show how to avoid the numerical error in the approximation of \eqref{eqn:generic:GLE} from deteriorating the accuracy and stability of the reduced switched system \eqref{eqn:sDAE:ROM}. In \Cref{sec:PMOR}, we propose a novel \MOR algorithm for \SLS and discuss its advantages with respect to the method of \cite{PonGB20}. Finally, we present numerical examples in \Cref{sec:examples}.

\subsection{Notation}
The symbol $\fI_l$, with $l\in\N$, denotes the identity matrix of size $\stateDim$. For $\fA\in\R^{\stateDim\times\stateDim}$, we write $\fA\succeq \zeroMat$ or $\fA \preceq \zeroMat$ if $\fA$ is symmetric positive or negative semidefinite, respectively. If $\fA$ has only real eigenvalues, then we denote the $i$-th largest eigenvalue (singular value) of $\fA$ by $\lambda_i(\fA)$ ($\sigma_i(\fA)$) and use $\kappa(\fA)=\tfrac{\sigma_1(\fA)}{\sigma_\stateDim(\fA)}$ to denote the spectral condition number of $\fA$. Finally, $\Lambda(\fA)$ indicates the spectrum of $\fA$.

\section{Preliminaries} \label{sec:prliminaries}
In this section, we recall the definitions of the reachable and observable sets for \SLS and then introduce the balancing-based \MOR algorithm for \SLS, which was first presented in \cite{PonGB20}.
\subsection{Reachability and Observability}
\label{subsec:reachObserv}
We start by introducing the concepts of reachable and observable sets for a continuous time-switching linear system. Let $\boldsymbol{\phi}( t, t_0, \stx_0, \inp, \switch)$ denote the state trajectory at time $t$ of the switched system \eqref{eqn:sDAE} starting from $\stx(t_0) = \stx_0$ with input $\inp$ and switching path $\switch\in\calS$ with $\calS$ given in \eqref{eqn:suitableSwitchingSignals}. 

\begin{definition}[Reachable and observable sets for \SLS]
	Let $\switch\in\calS$ be a given switching path. A state $\stx\in\R^{\stateDim}$ is called
	\begin{enumerate}
		\item \emph{reachable via $\switch$} if there exists a time instant $t_{\mathrm{f}} > t_0$, and an input $\inp\colon [t_0,t_f] \rightarrow \R^{\inpDim}$, such that $\boldsymbol{\phi}(t_f,t_0,0,\inp,\switch) = \stx$;
		\item \emph{unobservable via $\switch$}, if there exists an input $\inp$, such 
		\begin{equation*}
			\fC_{\switch}\boldsymbol{\phi}(t,t_0,\stx,\inp,\switch) = \fC_\switch \boldsymbol{\phi}(t, t_0,\zeroVec,\inp,\switch)\quad \text{for all } t \ge t_0.
		\end{equation*}
	\end{enumerate}
	The \emph{reachable} and \emph{unobservable sets} via $\switch$, denoted by $\calR_\switch$ and $\calU\calO_\switch$, respectively, are the sets of states that are reachable and unobservable via $\switch$, respectively. We define the \emph{observable set} via $\switch$ of \eqref{eqn:sDAE}, denoted by $\calO_\switch$, as $\calO_\switch \vcentcolon= (\calU \calO_\switch)^{\perp}$ (note that this definition is not unique since any complement set, not just the orthogonal one, would be suitable). The set of reachable states $\calR$ and the set of observable states $\calO$ of \eqref{eqn:sDAE} can be defined as
	\begin{align}
		\label{eqn:reachableObservableSet}
		\calR \vcentcolon= \bigcup_{\switch\in\calS}\calR_{\switch} \qquad\text{and}\qquad
		\calO \vcentcolon= \bigcup_{\switch\in\calS}\calO_{\switch}.
	\end{align}
\end{definition}
Reachability and observability are crucial in balancing-based \MOR. Indeed, the fundamental idea of this type of projection-based \MOR is to remove from the systems the states that are difficult to reach and/or difficult to observe; see \cite[Cha.~$7$]{Ant05} for further details.


\subsection{Model reduction algorithm: stability and accuracy conditions}
\label{sec:MORexplained}
Recall that the \MOR algorithm from \cite{PonGB20} obtains the projection matrices $\fV$ and $\fW$ by solving the \GLE~\eqref{eqn:GLE}. The symmetric positive-semidefinite solutions $\fcP,\fcQ\in\R^{\stateDim\times\stateDim}$ are known as the Gramians for~\eqref{eqn:sDAE} and they satisfy the relation
\begin{equation}\label{eqn:reach:observ:gramians}
\calR\;=\;\text{range}(\fcP)\quad\text{and}\quad\calO\;=\;\text{range}(\fcQ),
\end{equation}
with $\calR$ and $\calO$ as defined in \eqref{eqn:reachableObservableSet}; see \cite[Thm.~3]{PonGB20}. Let us emphasize that, at this stage, the choice of $\fA=\fA_1$ in \eqref{eqn:GLE} is somewhat arbitrary, and other choices of $\fA$ in the set $\{\fA_1,\ldots,\fA_M\}$ are also possible, provided that the $\fN_j$ matrices are then defined accordingly.

Next, compute (possibly singular) Cholesky factors $\fcP=\fS\fS^\T$ and $\fcQ=\fR\fR^{\T}$. Form the product $\fS^{\T}\fR$ and perform its singular value decomposition
\begin{equation}
	\label{eqn:Hank:mat}
	\fS^\T\fR = \left[\fU_1, \fU_2\right]\begin{bmatrix}
			\mathbf{\Sigma}_1&\zeroVec\\
			\zeroVec&\mathbf{\Sigma}_2
	\end{bmatrix}[\fV_1,\fV_2]^\T.
\end{equation}
Then, compute the projection matrices $\fV$ and $\fW$ via
\begin{equation}
	\label{eq:PGproj:mat}
	\fV = \fS\fU_1\mathbf{\Sigma}^{-1/2}_1 \quad\text{and}\quad
	\fW = \fR\fV_1\mathbf{\Sigma}^{-1/2}_1.
\end{equation}
This procedure is called square-root balanced truncation (see \cite[Sec.~7.3]{Ant05}). 

To comment on the stability of the \ROM derived through the procedure just detailed, we first need to briefly discuss the stability of the \SLS in \eqref{eqn:sDAE}; for a comprehensive overview on the topic, see \cite[Part~\MakeUppercase{\romannumeral 2}]{Lib03}. The first observation is that assuming each subsystem to be asymptotically stable, i.e., that the matrices $\fA_j$ are Hurwitz for all $j\in\switchingSet$, does not guarantee that the \SLS~\eqref{eqn:sDAE} will be stable for any $\switch\in\calS$. 
\begin{definition}\label{def:quad:stab}
	The switched linear system~\eqref{eqn:sDAE} is called \emph{quadratically stable} if there exists a positive definite matrix $\fP$ such that $\fA_j^{\T}\fP+\fP\fA_j \prec \zeroMat$, for all $j\in\switchingSet$.
\end{definition}
Quadratic stability plays a crucial role, as it offers a sufficient condition for ensuring the exponential stability of any \SLS characterized by stable modes, regardless of the chosen switching signals $\switch\in\calS$, where $\calS$ is specified in \eqref{eqn:suitableSwitchingSignals}. It is important to note that quadratic stability is merely a sufficient condition; a \SLS can remain stable for every $\switch\in\calS$ without necessarily being quadratically stable; see the discussion in \cite[Sec.~2.1.4]{Lib03}. The following lemma provides a crucial result for \MOR using the projection matrices in \eqref{eq:PGproj:mat}.
\begin{lemma}[{\!\!\cite[Lem.~12]{PetWL13}}]
	\label{lemma:qs:MOR}
	If the matrices $\fcP,\fcQ$ used for the construction of the projection matrices $\fV,\fW$ in~\eqref{eq:PGproj:mat} satisfy either 
	\begin{subequations}
		\label{eqn:LMI}
	\begin{equation}
		\label{eqn:LMI:P}
		\fA_j \fcP+\fcP\fA_j^{\T}+\fB_j\fB_j^\T \prec \zeroVec \quad \text{for all } j\in\switchingSet,
	\end{equation}
	or
	\begin{equation}
		\label{eqn:LMI:Q}
		\fA_j^\T \fcQ + \fcQ\fA_j+\fC_j^\T\fC_j \prec \zeroVec \quad\text{for all } j\in\switchingSet,
	\end{equation}
	\end{subequations}
	then the reduced system~\eqref{eqn:sDAE:ROM} obtained via $\fV$ and $\fW$ is quadratically stable. 
\end{lemma}
Finally, we recall a key result that offers an error estimate suitable for constructing a \ROM using balancing-based \MOR, ensuring that the \ROM output meets any prescribed accuracy requirements.
\begin{theorem}[{\!\!\cite[Thm.~6]{PetWL13}}]
	\label{thm:BTerrorBound}
	Let $\mathbf{\Sigma}_1$~in \eqref{eqn:Hank:mat} be of size $\stateDimRed\times \stateDimRed$, $\Lambda(\mathbf{\Sigma}_1)\cap\Lambda(\mathbf{\Sigma}_2) = \emptyset$, and let the distinct singular values of $\mathbf{\Sigma}_2$~be $\sigma_k$, with multiplicities $m_k $, for  $k= \stateDimRed+1, \ldots ,\tilde \stateDim$.  It follows $\tilde \stateDim \le \stateDim$ and $\stateDim=\stateDimRed+\sum_{k=\stateDimRed+1}^{\tilde \stateDim}m_i$. Suppose $\fcP$ and $\fcQ$ satisfy the \LMIs
	\begin{equation}
		\label{eqn:LMI:ERR}
		\fA_j \fcP+\fcP\fA_j^{\T}+\fB_j\fB_j^\T \preceq \zeroVec \quad\text{and}\quad\fA_j^\T \fcQ + \fcQ\fA_j+\fC_j^\T\fC_j \preceq \zeroVec \quad \text{for all } j\in\switchingSet.
	\end{equation}
	Then, for every switching signal $\switch\in\calS$, the output error between the \FOM~\eqref{eqn:sDAE} and the \ROM~\eqref{eqn:sDAE:ROM}, obtained through the projection matrices \eqref{eq:PGproj:mat}, is bounded by twice the sum of the neglected singular values, multiplicities not counted, times the $L_2$ norm of the input function, i.e.,
	\begin{equation}
	\label{eq:err:est}
	\|\out-\hat{\out}\|_{L_2} \leq 2 \left(\sum_{k=\stateDimRed+1}^{\tilde \stateDim}\sigma_{k}\right) \|\inp\|_{L_2}.
\end{equation}
\end{theorem}
\begin{proof}
The proof is derived from \cite[Thm.~6]{PetWL13}. The unique aspect here is that the original authors did not specifically differentiate between simple and non-simple neglected singular values. Nevertheless, tailoring their proof to encompass this distinction is achieved by reiterating their steps while addressing non-simple singular values. This approach is further elucidated in the subsequent proof of \Cref{teo5}, so those steps will not be reiterated here.
\end{proof}

We conclude that for the quadratic stability of the reduced \SLS, compliance with just one set of \LMIs is necessary, specifically either \eqref{eqn:LMI:P} or \eqref{eqn:LMI:Q}. To ensure error guarantees in the reduced approximation, the sets of \LMIs presented in condition \eqref{eqn:LMI:ERR} must also be fulfilled. These \LMIs are generally less stringent than those for quadratic stability. Consequently, quadratic stability does not inherently ensure that \eqref{eq:err:est} is satisfied; conversely, satisfying~\eqref{eq:err:est} does not inherently ensure the quadratic stability of \eqref{eqn:sDAE:ROM}.

\section{About numerics for large-scale GLEs}
\label{sec:numericalDetails}

For the large-scale setting, it is crucial to exploit sparsity when available and to use state-of-the-art methods in every step of the \MOR procedure presented in \Cref{sec:MORexplained}. Besides computational efficiency, we further aim for error certification and thus have to balance the errors of the different approximations. In more detail, we discuss the algorithm used to approximate the solutions of \eqref{eqn:GLE} and derive an error estimate to certify the quality of the approximation in \Cref{subsec:app:GLE}. Afterward, in \Cref{subsec:numerics:sec3} we examine the numerical efficiency of computing these error bounds, explicitly identifying scenarios in which we can attain both efficiency and error control, as well as scenarios in which only one of these properties can be guaranteed. We conclude with \Cref{subsec:red:err} where we show how the provided error estimate can be used to prevent instabilities and deterioration of accuracy in the \ROM for a certain class of \SLS.

\subsection{Solvability and stationary algorithm}
\label{subsec:stationaryAlgo}
Let us write the \GLE in the generic form
\begin{equation}\label{eqn:generic:GLE}
\fA\fX+\fX\fA^{\T}+\sum_{i=1}^{M}\left(\fN_j\fX\fN_j^{\T}\right)+\fB\fB^{\T},
\end{equation}
where $\fA$ is Huwritx and $\fB\in\R^{\stateDim\times\tilde\inpDim}$, with $\tilde \inpDim\ll\stateDim$.
The existence and uniqueness of a solution of~\eqref{eqn:generic:GLE} can be studied via Kronecker
products, i.e., by introducing the matrices
\begin{align}\label{eqn:GLE:oper}
       \lyapOper \vcentcolon= \fI_\stateDim\otimes\fA +\fA\otimes\fI_\stateDim \in\R^{\stateDim^2\times\stateDim^2} \qquad\text{and}\qquad
       \gleOper \vcentcolon= \sum_{j=1}^{M}\fN_j\otimes\fN_j \in\R^{\stateDim^2\times\stateDim^2},
\end{align}
such that the vectorized form of~\eqref{eqn:generic:GLE} is given by $(\lyapOper+\gleOper)\vec2(\fX)=-\vec2(\fB\fB^\T)$.
Hence, the \GLE~\eqref{eqn:generic:GLE} is uniquely solvable if and only if $\lyapOper+\gleOper$ is nonsingular. 

Instead of solving the linear system directly, we employ the stationary iteration algorithm from \cite[Alg.~2.1]{ShaSS16} to solve the \GLE~\eqref{eqn:generic:GLE} numerically; see \Cref{alg:statIterGLE}. 
\begin{algorithm}[ht]
	\caption{Stationary iterations for the \GLEs~\eqref{eqn:generic:GLE}}
	\label{alg:statIterGLE}
	\begin{minipage}{\linewidth}
		\textbf{Input:} Matrices $\fA$, $\fB$,  $\fN_j$ for $j=1,\ldots,M$\\
		\textbf{Output:} $\tilde{\fZ}$ with $\tilde \fX=\tilde{\fZ}\tilde{\fZ}^\T$ 
	\end{minipage}
	\begin{algorithmic}[1]
		\State Set ${\fB}_1 \vcentcolon=\fB$
		\State Approximately solve $\fA\fX+\fX\fA^\T+{\fB}_1{\fB}_1^\T=\zeroMat$ for $\fX_1=\fZ_1\fZ_1^\T$\label{line2}
		\For{$k=2,3,\ldots$}
			\State Set ${\fB}_k \vcentcolon= [\fN_1\fZ_{k-1},\ldots,\fN_M\fZ_{k-1},{\fB}_1]$
			\State Approximately solve \label{A2-l5}
				\begin{equation}
					\label{alg:statIterGLE:lyap}
					\fA\fX+\fX\fA^\T+{\fB}_k{\fB}^\T_k = \zeroMat
				\end{equation}
				\phantom{\textbf{for}}for $\fX_k=\fZ_k\fZ_k^\T$
			\State \textbf{if} {sufficiently accurate} \textbf{then} stop
		\EndFor 
	\end{algorithmic}
\end{algorithm}
It consists of a fixed-point iteration in which we need to solve a classical Lyapunov equation in each iteration. To ensure that \Cref{alg:statIterGLE} is feasible and convergent, we make the following assumption; see, for instance, \cite{JarMPR18}.

\begin{assumption}
	\label{ass:statIterGLE}
	The matrix $\fA\in\R^{\stateDim\times\stateDim}$ in~\eqref{eqn:generic:GLE} is Hurwitz and $\|\lyapOper^{-1}\gleOper\|_2<1$.
\end{assumption}

\begin{remark}
	\label{rem:rescaling}
	On first glance, the assumption $\|\calL^{-1}\gleOper\|_2<1$ sounds rather restrictive. Nevertheless, in the context of \MOR, we can always proceed by defining the rescaled matrix
	\begin{equation*}
		\tilde{\gleOper} = \frac{1}{\|\calL^{-1}\gleOper\|_2+\delta}\gleOper,
	\end{equation*}
	for some $\delta>0$ and observing that rescaling $\gleOper$ does not affect the space spanned by the columns of the solution matrix $\fX$; see \cite{ConI05,BenD11}.  Since the columns'space stays the same, the reachability and observability sets of the \SLS are also correctly captured by the column space of the solution of the rescaled system; see relation \eqref{eqn:reach:observ:gramians}. We postpone the discussion of the estimation of $\|\calL^{-1}\gleOper\|_2$ in \Cref{subsub:gamma:sc}.
\end{remark}

\subsection{Error certification}
\label{subsec:app:GLE}

An efficient implementation of \Cref{alg:statIterGLE} involves solving the following numerical issues:
\begin{enumerate}
	\item At each iteration, the Lyapunov equation requires an efficient solver.
	\item An opportune stopping criterion is required for the stationary iterations.
\end{enumerate}
In the following subsections, we describe how to address these points.

\subsubsection{Approximate solution of Lyapunov equations involving large sparse matrices}
\label{subsubsec:Lyapunov}

Consider \eqref{alg:statIterGLE:lyap} with a possibly large and sparse matrix $\fA$. The literature on the case of large and sparse matrices is vast; see, for example, the survey \cite{Sim16} on methods for matrix equations and the references therein. Almost all methods in the large-scale setting have in common that they construct a sequence of approximations that converge to the true solution. Thus, a suitable stopping criterion is required to balance the approximation error in the calculation of the solution of the Lyapunov equation and the solution of the stationary iteration from \Cref{alg:statIterGLE}. To this end, assume that $\fX_{\indKS}\in\Spsd{\stateDim}$ is the approximation of the solution $\fX\in\Spsd{\stateDim}$ of the Lyapunov equation~\eqref{alg:statIterGLE:lyap} in the $\indKS$th step (and $k$th iteration of \Cref{alg:statIterGLE}) and let
\begin{equation}
	\label{eqn:lyap:residual}
	\fR_\indKS \vcentcolon= \fA\fX_\indKS + \fX_\indKS\fA^\T + {\fB}_k{\fB}_k^\T
\end{equation}
denote the corresponding residual. We obtain the following error-residual relation.

\begin{proposition}[Error bound for approximate solution of Lyapunov equations]
	\label{prop1}
	Consider the Lyapunov equation~\eqref{alg:statIterGLE:lyap} with $\fA$ Hurwitz, and assume that $\fX_\indKS\in\Spsd{\stateDim}$ is an approximation of the unique solution $\fX\in\Spsd{\stateDim}$ of~\eqref{alg:statIterGLE:lyap}. Let $\fE_\indKS \vcentcolon= \fX_\indKS-\fX$. Then
	\begin{equation}
		\label{eqn:lyap:errorResidual}
		\|\fE_\indKS\|_2 \leq \frac{\|\fR_\indKS\|_{\Frob}}{\sigma_{\min}(\lyapOper)},
	\end{equation}
	where $\sigma_{\min}(\fA)$ denotes the smallest singular value of the matrix $\fA$.
\end{proposition}

\begin{proof}
	Substituting~\eqref{alg:statIterGLE:lyap} into the residual equation~\eqref{eqn:lyap:residual} yields the error-residual relation
	\begin{align}
		\label{eqn:errResLyap}
		\fR_\indKS = \fA \fE_\indKS + \fE_\indKS\fA^\T,
	\end{align}
	which, together with the symmetry of $\fR_\indKS$, yields that $\fE_\indKS$ is symmetric; see, e.g., \cite[Cha.~12.3, Thm.~3]{LanT85}. Using Kronecker algebra and the notation introduced in \eqref{eqn:GLE:oper}, we can reformulate~\eqref{eqn:errResLyap} in vectorial form as $\vec2{(\fR_\indKS)} = \calL \vec2{(\fE_\indKS)}$.
	Using $\|\fE\|_2\le\|\fE\|_{\Frob}$ yields
		\begin{align*}
	   \|\fE_\indKS\|_2 \;&\leq \;\|\fE_\indKS\|_{\Frob} \;= \;\|\vec2(\fE_\indKS)\|_2 \;=\; \|\calL^{-1}\vec2{(\fR_\indKS)}\|_2 \;\leq\; \|\calL^{-1}\|_2\|\vec2{(\fR_\indKS)}\|_2\\
	   &=\; \frac{\|\fR_\indKS\|_{\Frob}}{\sigma_{\min}(\calL)}.\qedhere
	\end{align*}
\end{proof}

\subsubsection{Stopping criteria for \texorpdfstring{\Cref{alg:statIterGLE}}{Algorithm 2}}

Consider the approximate solution $\fX_k$ for the \GLE~\eqref{eqn:generic:GLE} provided by \Cref{alg:statIterGLE} at step $k$, and assume for the moment that the Lyapunov equation \eqref{alg:statIterGLE:lyap} is solved exactly. To understand whether such a solution is accurate enough for our requirements, we must estimate the norm of the approximation error.

\begin{proposition}
	\label{prop2}
	Suppose \Cref{ass:statIterGLE} is satisfied and consider the error $\fE_k \vcentcolon=\fX-\fX_k$, where $\fX$ is the unique solution of the \GLE~\eqref{eqn:generic:GLE} and $\fX_k = \fZ_k\fZ_k^\T$ is the unique solution of the Lyapunov equation~\eqref{alg:statIterGLE:lyap} in iteration $k$ of \Cref{alg:statIterGLE}. Then \Cref{alg:statIterGLE} converges and
	\begin{equation}\label{eqn:stopping:GLE}
		\|\fE_k\|_2 \leq \gamma \|\fX_k-\fX_{k-1}\|_{\Frob},
	\end{equation}
	with
	\begin{equation}
		\label{eqn:GLEerrorBoundConstant}
		\gamma \vcentcolon= \frac{\|\calL^{-1}\gleOper\|}{1-\|\calL^{-1}\gleOper\|}.
	\end{equation}
\end{proposition}

\begin{proof}
	Due to \Cref{alg:statIterGLE} we obtain	
	\begin{equation}\label{eq30}
		\fA\fX_k + \fX_k\fA^\T + \sum_{j=1}^{M}\left(\fN_j\fX_{k-1}\fN_j^\T\right) + {\fB}_{1}{\fB}_{1}^\T = \zeroMat
	\end{equation}
	for all $k\geq 2$. Then, from \eqref{eqn:generic:GLE}, we have ${\fB}_{1}{\fB}_{1}^\T=-\fA\fX-\fX\fA^\T-\sum_{j=1}^{M}\left(\fN_j\fX\fN_j^\T\right)$, and substituting into \eqref{eq30}, we get
	\begin{equation}\label{eq31}
		\fA\fE_k+\fE_k\fA^\T+\sum_{j=1}^{M}\left(\fN_j\fE_{k-1}\fN_j^\T\right)=\zeroMat.
	\end{equation}
	For convenience, we write \eqref{eq31} in its vectorial form $\calL\vec2(\fE_k) = -\gleOper \vec2(\fE_{k-1})$ with~$\calL$ and~$\gleOper$ defined in \eqref{eqn:GLE:oper}; thus, we have
	\begin{equation}\label{eq32}
		\vec2(\fE_k)=-\calL^{-1}\gleOper\vec2(\fE_{k-1}).
	\end{equation}
	Relation \eqref{eq32} directly implies that, from a numerical point of view, \Cref{alg:statIterGLE} converges, i.e., $\lim_{k\rightarrow\infty}\|\fE_k\|_2=0$, if $\rho(\calL^{-1}\gleOper)<1$, where $\rho(\calL^{-1}\gleOper)$ denotes the spectral radius of $\calL^{-1}\gleOper$; see also \cite{Dam08,ShaSS16}. Observe that $\fE_{k-1}=\fE_{k}+(\fX_{k}-\fX_{k-1})$, therefore, substituting in \eqref{eq32} and considering the norm on both sides, we get
	\begin{align*}
		\begin{aligned}
			\|\fE_k\|_{\Frob} = \|\vec2(\fE_k)\|_2 &= \|\calL^{-1}\gleOper\vec2\left(\fE_k+\fX_{k}-\fX_{k-1}\right)\|_2\\
			&\le \|\calL^{-1}\gleOper\|_2{\Big(\|\fE_k\|_{\Frob}+\|\fX_{k}-\fX_{k-1}\|_{\Frob}\Big)}.
		\end{aligned}
	\end{align*}
	By \Cref{ass:statIterGLE}, we have $\|\calL^{-1}\gleOper\|<1$, and hence
	\begin{equation*}
		\|\fE_k\|_2 \le \|\fE_k\|_\Frob \le \frac{\|\calL^{-1}\gleOper\|_2}{1-\|\calL^{-1}\gleOper\|_2}\|\fX_{k}-\fX_{k-1}\|_\Frob.\qedhere
	\end{equation*}
\end{proof}
Lastly, we have to compute $\|\fX_j-\fX_{j-1}\|_\Frob$, which will require $\calO(\stateDim^2)$ floating point operations. To overcome this issue, we employ the following result.
\begin{proposition}
	\label{prop3}
	Let $\fX_{k}=\fZ_k\fZ^\T_{k}$ with $\fZ_k\in \R^{\stateDim\times \stateDim_k}$ and $\fX_{k-1}=\fZ_{k-1}\fZ^\T_{k-1}$ with $\fZ_{k-1}\in\R^{\stateDim\times \stateDim_{k-1}}$ denote the respective solutions at iterations $k$ and $k-1$ of \Cref{alg:statIterGLE}. Without loss of generality, assume $\stateDim_k = \stateDim_{k-1}$\footnote{If this does not hold, it is enough to fill the $|\stateDim_k-\stateDim_{k-1}|$ missing columns of the smaller matrix with zero entries.}. Let 
	\begin{equation*}
		\fX_k=\fU_k\fSigma_k\fU_k^{\T},\quad\quad\fX_{k-1}=\fU_{k-1}\fSigma_{k-1}\fU_{k-1}^{\T},
	\end{equation*}
	denote the truncated eigenvalue decomposition for $\fX_k$ and $\fX_{k-1}$, respectively, and assume
	\begin{equation*}
		\fZ_k =\fU_k\fSigma^{\frac{1}{2}}_k,\quad\quad \fZ_{k-1}=\fU_{k-1}\fSigma^{\frac{1}{2}}_{k-1}.
	\end{equation*}
	Then
	\begin{equation}
		\label{eq:46}
		\|\fX_k-\fX_{k-1}\|_\Frob \le \|(\fZ_k-\fZ_{k-1})\fSigma_k^{\frac{1}{2}}\|_\Frob+\|(\fZ_k-\fZ_{k-1})\fSigma_{k-1}^{\frac{1}{2}}\|_\Frob.
	\end{equation}
\end{proposition}

\begin{proof}
	By making use of the triangular inequality, it directly follows
	\begin{align*}
					\|\fX_k-\fX_{k-1}\|_\Frob &= \|\fZ_k\fZ_k^\T-\fZ_k\fZ_{k-1}^\T+\fZ_k\fZ_{k-1}^\T -\fZ_{k-1}\fZ_{k-1}^\T \|_\Frob\\
					&\le \|\fU_k\fSigma^{\frac{1}{2}}_k (\fZ_k^\T-\fZ_{k-1}^\T)\|_\Frob+\|(\fZ_k-\fZ_{k-1})\fSigma^{\frac{1}{2}}_{k-1}\fU^{\T}_{k-1}\|_\Frob\\
					&= \|(\fZ_k-\fZ_{k-1})\fSigma_k^{\frac{1}{2}}\|_\Frob+\|(\fZ_k-\fZ_{k-1})\fSigma_{k-1}^{\frac{1}{2}}\|_\Frob,
	\end{align*}
	where, for the last equality, we used the fact that $\fU_k$ and $\fU_{k-1}$ have orthonormal columns and that $\|\fSigma_k^{\frac{1}{2}}(\fZ_k-\fZ_{k-1})\|_\Frob=\|\fSigma_k^{\frac{1}{2}}(\fZ_k-\fZ_{k-1})^{\T}\|_\Frob$.
\end{proof}
	
Replacing $\|\fX_k-\fX_{k-1}\|_\Frob$ with the right-hand side of \eqref{eq:46} has the advantage that we only need $\calO(\stateDim)$ floating point operations for its evaluation (assuming $\stateDim_k,\stateDim_{k-1}\ll \stateDim$). With these preparations, we are now ready to formulate our main result for the approximate solution of the \GLE~\eqref{eqn:generic:GLE}.
	
\begin{theorem}[Error bound for the approximate solution of the \GLE]
	\label{teo1}
	Suppose \Cref{ass:statIterGLE} is satisfied. Let $\fX=\fZ\fZ^\T$ with $\fZ\in\R^{\stateDim\times\stateDim}$ be the unique solution of the \GLE \eqref{eqn:generic:GLE} and $\tilde{\fX}_k=\tilde{\fZ}_k\tilde{\fZ}_k^\T$ with $\tilde{\fZ}_k\in\R^{\stateDim\times\stateDim_k}$ be the approximate solution of the \GLE computed at iteration $k$ of \Cref{alg:statIterGLE} applying the an iterative method to solve the Lyapunov equation in~\eqref{alg:statIterGLE:lyap}. Then
	\begin{align}
		\label{eq:50}
		\|\fX-\tilde{\fX}_k\|_2 \le \gamma \Delta_{k,k-1}+ \frac{(1+\gamma)\|\fR_k\|_{\Frob}+\gamma\|\fR_{k-1}\|_{\Frob}}{\sigma_{\min}(\lyapOper)},
	\end{align}
	where $\fR_k$ and $\fR_{k-1}$ are the residuals corresponding to the approximation of the Lyapunov equation~\eqref{alg:statIterGLE:lyap} at the $k$th and $(k-1)$th iterations of \Cref{alg:statIterGLE} due to the iterative method for \eqref{alg:statIterGLE:lyap}, respectively, $\gamma$ is as defined in \eqref{eqn:GLEerrorBoundConstant}, and
	\begin{equation}
	  \Delta_{k,k-1} \vcentcolon= \|(\tilde \fZ_k-\tilde \fZ_{k-1})\tilde \fSigma_k^{\frac{1}{2}}\|_\Frob+\|(\tilde\fZ_k-\tilde\fZ_{k-1})\tilde\fSigma_{k-1}^{\frac{1}{2}}\|_\Frob.
	\end{equation}
\end{theorem}

\begin{proof}
	Let $\fX_k$ be the exact solution of the Lyapunov equation~\eqref{alg:statIterGLE:lyap} computed at step $k$ of \Cref{alg:statIterGLE}. Using \Cref{prop1,prop2} we obtain 
	\begin{align}
		\|\fX-\tilde{\fX}_k\|_2 &\le \|\fX-\fX_k\|_2 + \|\fX_k-\tilde{\fX}_k\|_2 \le \gamma\|\fX_k-\fX_{k-1}\|_\Frob + \frac{\|\fR_k\|_{\Frob}}{\sigma_{\min}(\lyapOper)}.
	\end{align}
	Furthermore, the triangle inequality and \Cref{prop1,prop3} imply
	\begin{align*}
		\|\fX_k-\fX_{k-1}\|_\Frob &\le \|\fX_k-\tilde{\fX}_k\|_\Frob+\|\tilde{\fX}_k-\tilde{\fX}_{k-1}\|_\Frob+\| \tilde{\fX}_{k-1}-\fX_{k-1}  \|_\Frob \nonumber\\
				&\le \frac{\|\fR_k\|_{\Frob}}{\sigma_{\min}(\lyapOper)}+\Delta_{k,k-1}+\frac{\|\fR_{k-1}\|_{\Frob}}{\sigma_{\min}(\lyapOper)},
	\end{align*}
	yielding \eqref{eq:50}.
\end{proof}
\subsection{Remarks on the numerical efficiency of rigorous stopping criteria}\label{subsec:numerics:sec3}
The practical use of \eqref{eq:50} as a stopping criterion for \Cref{alg:statIterGLE} requires clarification on how $\|\fR_k\|_{\Frob}$, $\sigma_{\min}(\lyapOper)$, and $\gamma$ are evaluated. Ideally, when addressing sparse problems of dimension $\stateDim$, the objective is to execute \Cref{alg:statIterGLE} using $\calO(\stateDim)$ floating-point operations, including the computations required to evaluate stopping criteria that rigorously certify the accuracy of the resulting approximation. In this section, we show that such a complexity bound does not appear to be achievable, in general, when also accounting for rigorous error control. 

\subsubsection{Krylov subspace methods for the approximation of Lypunov equations}

 In this section, we detail the following: ($i$) the choice of the iterative method to determine $\fX_\indKS$, and ($ii$) how to efficiently compute the residual norm $\|\fR_\indKS\|_{\Frob}$. Since these questions are closely tied together, we first briefly describe the method for computing $\fX_\indKS$. As proposed in \cite{Sim16}, we project the Lyapunov equation~\eqref{eqn:lyap:residual} onto a smaller subspace and then use standard solvers for the small equation, thus obtaining the coordinates of the approximation in the subspace. We refer to \cite[Sec.~5.1]{Sim16} for a review of the solution of small-scale Lyapunov equations. 

More precisely, consider the matrix $\fV_\indKS\in\R^{\stateDim\times \stateDim_\indKS}$ whose $\stateDim_\indKS$ columns form an orthonormal basis of the appropriate subspace $\calV_\indKS$ that is to be used for the projection. Let ${\fX}_\indKS=\fV_\indKS\fY_\indKS\fV_\indKS^\T$ be a low-rank approximation to the solution of \eqref{eqn:lyap:residual}  with $\fY_\indKS\in\R^{\stateDim_\indKS \times \stateDim_\indKS}$. Imposing the Galerkin condition $\fV_\indKS^\T\fR_\indKS\fV_\indKS = \zeroMat$ on the residual, we obtain the lower-order projected Lyapunov equation
\begin{equation}
	\label{eqn:lyap:projected}
	\fA_\indKS\fY_\indKS + \fY_\indKS\fA_\indKS^\T + {\fB}_\indKS{\fB}_\indKS^\T = \zeroMat,
\end{equation}
with $\fA_\indKS \vcentcolon= \fV_\indKS^\T\fA\fV_\indKS$ and ${\fB}_\indKS \vcentcolon= \fV^\T_\indKS{\fB}_k$. For $\calV_\indKS$, one can choose, for instance, the proposal from \cite{Saa89} and use the block Krylov subspace
\begin{equation}\label{eq:krol:def}
	\calK_\indKS(\fA,{\fB}_k) = \spann \{{\fB}_k,\fA{\fB}_k,\fA^2{\fB}_k,\ldots,\fA^{\indKS-1}{\fB}_k\}.
\end{equation}
For more sophisticated techniques, such as rational Krylov subspaces, we refer to \cite{Sim07,Sim16}.

Besides computational efficiency, the low-rank approximation can also be exploited for an efficient evaluation of the Frobenius norm of the residual.  Indeed, using the Krylov subspace \eqref{eq:krol:def}, one can show $\|\fR_\indKS\| = \|\fV_{\indKS+1}^\T \fR_{\indKS} \fV_{\indKS+1}\|$ for the spectral and the Frobenius norm; see \cite[Prop.~3.3]{Sim07}. This result is particularly convenient since, for $\fV_{\indKS+1}\vcentcolon =[\fV_\indKS,\;\tilde{\fV}]$, we have
\begin{align}\label{eq:61}
	\fV^\T_{\indKS+1}\fR_{\indKS}\fV_{\indKS+1}=\begin{bmatrix}
		\fV_{\indKS}^\T\fR_{\indKS}\fV_{\indKS} & \fV_{\indKS}^\T\fR_{\indKS}\tilde{\fV}\\
\tilde{\fV}^\T \fR_{\indKS}\fV_{\indKS}	& \tilde{\fV}^\T\fR_{\indKS}\tilde{\fV}
	\end{bmatrix}=\begin{bmatrix}
	\zeroMat &\fY_\indKS\fV_{\indKS}^\T\fA^\T	\tilde{\fV} \\
	 \tilde{\fV}^\T\fA\fV_{\indKS}\fY_\indKS&\zeroMat
	\end{bmatrix}
\end{align}
where we used the expression \eqref{eqn:lyap:residual}, the fact that ${\fX}_\indKS=\fV_\indKS\fY_\indKS\fV_\indKS^\T$, and that $\tilde{\fV}$ is, by construction, orthogonal to $\fV_\indKS$ and ${\fB}_k$. From \eqref{eq:61}, it directly follows that
\begin{equation}\label{eqn:eff:res}
	\|\fR_{\indKS}\|_{\Frob} = \| \fV^\T_{\indKS+1}\fR_{\indKS}\fV_{\indKS+1}  \|_\Frob =\;2\| \tilde{\fV}^\T\fA\fV_{\indKS}\fY_\indKS \|_\Frob
\end{equation}
and the right hand side of \eqref{eqn:eff:res} can be used to evaluate $\|\fR_{\ell}\|_\Frob$ in a computationally efficient way without the need to ever form the large-scale matrix $\fR_{\ell}$. Finally, assuming the subspace procedure for iteration $k$ of \Cref{alg:statIterGLE} has converged to $\fY_k$, to use \Cref{prop2}, we need the truncated eigenvalue decomposition related to $\fX_k$ to form $\fZ_k$. This is efficiently obtained if $\fY_k$ results in a small dimension, by computing the eigenvalue decomposition of the corresponding matrix $\fY_k$, i.e., $\fY_k=\fU_{\fY}\fSigma_{\fY}\fU_{\fY}^{\T}$, and then we can set $\fU_k=\fV_k\fU_{\fY}$ and $\fSigma_k=\fSigma_{\fY}$.

\subsubsection{On the smallest singular value of the Lyapunov operator}\label{subsubsec:eff:sig:min}
We detail the computational complexity associated with evaluating $\sigma_{\min}(\lyapOper)$, with $\lyapOper$ as defined in \eqref{eqn:GLE:oper}. To this end, we proceed from specific instances to increasingly general settings. Please note that none of the subsequent bullet points introduce novel results. Rather, their purpose is to provide an explicit and exhaustive characterization that depends on the properties of the matrix $\fA$.
\begin{itemize}
    \item If $\fA$ is symmetric and negative definite (thus Hurwitz), then $\sigma_{\min}(\lyapOper)=2\lambda_{\min}(-\fA)$. Thus, the computational cost for sparse matrices scales as $\calO(\stateDim)$. 
    \item If $\fA$ is normal, i.e., $\fA_j\fA^{*}_j=\fA_j^{*}\fA_j$, and Hurwitz, then 
    \begin{equation*}
        \sigma_{\min}(\lyapOper)\;=\;2\lambda_{\min}\left(-\frac{\fA+\fA^{*}}{2}\right),
    \end{equation*}
 i.e., the smallest singular value of the Lyapunov operator is equal to two times the smallest eigenvalue of the Hermitian part of $-\fA$. The computation of this eigenvalue via iterative methods scales linearly with the state dimension $\stateDim$.
 \item If $-\fA$ is coercive, that is, if $-(\fA+\fA^{*})\succ 0$, or equivalently, if $\fA$ is a strictly dissipative matrix, then
 \begin{equation*}
        \sigma_{\min}(\lyapOper)\;\ge\;2\lambda_{\min}\left(-\frac{\fA+\fA^{*}}{2}\right).
    \end{equation*}
  \item For a generic matrix $\fA$, there are two options: 
  \begin{enumerate}
      \item computing $\sigma_{\min}(\lyapOper)$ using an iterative method tailored to sparse matrices, which, however, would entail $\calO(\stateDim^2)$ floating-point operations; see for instance \cite{SchS03}. We emphasize that computing $\sigma_{\min}(\lyapOper)$ remains more advantageous than solving \eqref{eqn:generic:GLE} by forming the $\stateDim^2$ linear systems, since the latter procedure would require explicitly storing $\fX$ (full square matrix of dimension $\stateDim$);
      \item alternative bounds on $\|\fE_{\indKS}\|_2$, that do not involve the computations of $\sigma_{\min}(\lyapOper)$, can be derived using pseudospectral knowledge of $\fA$; see \Cref{app:A}.  
\end{enumerate}
\end{itemize}
\subsubsection{About $\gamma$ and the stopping criteria}\label{subsub:gamma:sc}
For the computation or approximation of $\|\calL^{-1}\gleOper\|_2$, we observe
\begin{align}
	\label{eq33}
	\|\calL^{-1}\gleOper\|_2 &\leq 
	\|\calL^{-1}\|_2\|\gleOper\|_2 = \frac{1}{\sigma_{\min}(\lyapOper)}\Bigg\|\sum_{j=1}^{M}\fN_j\otimes\fN_j\Bigg\|_2\nonumber\\
	&\le \frac{1}{\sigma_{\min}(\lyapOper)}\sum_{j=1}^{M} \Big\|\fN_j\otimes\fN_j\Big\|_2 = \frac{1}{\sigma_{\min}(\lyapOper)}\sum_{j=1}^{M}\sigma_1(\fN_j)^2,
\end{align}
and \eqref{eq33} computational cost depends on the one for $\sigma_{\min}(\lyapOper)$. To ensure the convergence of \Cref{alg:statIterGLE}, for any $\delta>0$, we perform the scaling
\begin{equation}
	\label{eq:52}
	\tilde{\fN}_j \vcentcolon=\frac{\fN_j}{\sqrt{\beta+\delta}},\quad\text{with}\quad \beta \coloneqq \frac{1}{\sigma_{\min}(\lyapOper)}\sum_{j=1}^{M}\sigma_1(\fN_j)^2,\quad j=1,\ldots M;
\end{equation}
see also \Cref{rem:rescaling}. Hence, it is not restrictive to consider $\gamma\approx 1$. Indeed, it is sufficient to scale \eqref{eq:52} with $\delta = \beta$ to get $\gamma<1$, and thus we can get rid of~$\gamma$ in~\eqref{eqn:stopping:GLE}.\label{enum:1}

Considering $\gamma=1$, we can require
\begin{align*}
	\|\fR_k\|_\Frob \le \zeta_1\frac{\sigma_{\min}(\lyapOper)}{3}\tol, \quad\quad
\Delta_{k,k-1}\le 	\zeta_2\tol,\quad\text{with}\quad\zeta_1+\zeta_2=1
\end{align*}
during \Cref{alg:statIterGLE} to obtain $\|\fX-\fX_k\|_2\le \tol$. Our default choices are $\zeta_2=0.9$ and $\zeta_1=0.1$. Otherwise, a clever implementation consists of dynamically updating the exit tolerance for the approximation of \eqref{alg:statIterGLE:lyap}. For instance, $\|\fR_k\|_\Frob \le\tol_{\mathrm{Lyap},k} $ where
\begin{equation*}
\tol_{\mathrm{Lyap},k}\;\vcentcolon =\; \max\left\{\zeta_1\frac{\sigma_{\min}(\lyapOper)}{3}\tol,\min\left\{ \frac{\zeta_1}{\zeta_2}\Delta_{k-1,k-2},\tol_{\mathrm{Lyap},k-1}\right\}\right\}
\end{equation*}
This allows us to save some computational effort since the accuracy for the resolution of the Lyapunov equation is scaled based on the current level of accuracy reached by the stationary iteration algorithm.

\subsection{The \GLE error approximation and \MOR}
\label{subsec:red:err}
The error bound from \Cref{thm:BTerrorBound} relies on the  singular values defined in \eqref{eqn:Hank:mat} and on the \LMIs~\eqref{eqn:LMI}. Both depend on the solution of the \GLEs~\eqref{eqn:GLE}, which we can solve to a given tolerance using the results from the previous subsection. Nevertheless, \Cref{thm:BTerrorBound} relies on the exact Gramians, so it is fundamental to understand the impact of the numerical approximation error. 

\subsubsection{Preservation of \LMIs for the approximated solution}
\label{subsection:3.2.1}

After running algorithm \Cref{alg:statIterGLE} for k iterations, we obtain an approximate solution $\fX_k$ for the exact solution $\fX$. The central problem we address is whether this approximation preserves the stability of the \ROM and ensures that the error bound stated in \Cref{thm:BTerrorBound} remains valid. Specifically, even if the exact Gramians satisfy the \LMIs in \eqref{eqn:LMI:ERR} and one of the \LMIs in \eqref{eqn:LMI}, there is no guarantee that the computed approximations will do the same. Therefore, we aim to establish conditions under which, if the exact Gramian $\fX$ (either the reachability Gramian satisfying \eqref{eqn:LMI:P} or the controllability Gramian satisfying \eqref{eqn:LMI:Q}) fulfills the \LMIs, then our approximation $\fX_k$ also does. This would guarantee the quadratic stability of the \ROM and the applicability of \Cref{thm:BTerrorBound}.

With this in mind, let us introduce the matrix $\hat{\fX}$ that, under certain conditions on the~$\fA_j$ matrices, satisfies the same \LMIs as the exact solution~$\fX$. We construct $\hat{\fX}$ using the approximation $\tilde{\fX}$ given by \Cref{alg:statIterGLE} and leveraging \Cref{teo1}.

\begin{theorem}\label{teo4}
	Suppose that the exact solution $\fX$ of~\eqref{eqn:GLE:reach} satisfies the \LMI \eqref{eqn:LMI:P} and let $\tilde{\fX}$ denote its approximation given by \Cref{alg:statIterGLE}. Assume that the $\fA_j$ are strictly dissipative, i.e.,
	\begin{equation}\label{eqn:ass:ndeAj}
	\fA_j+\fA^{\T}_j\prec\zeroVec,\quad\text{ for all }j\in\switchingSet.
	\end{equation}
	Let $\perturbation\geq 0$ satisfy
	\begin{align}\label{eqn:ass:teo4}
		\begin{aligned}
		\perturbation \;\ge\; \max_{j\in\switchingSet} \frac{\sigma_1(\fA_j)}{\sigma_\stateDim(\fA_j)} \|\fX-\tilde \fX\|_2\;=\;\max_{j\in\switchingSet}\kappa(\fA_j)\|\fX-\tilde \fX\|_2,
		\end{aligned}
	\end{align}
	where $\kappa(\fA_j)$ denotes the spectral condition number of $\fA_j$. Then, the matrix
	\begin{equation}\label{eqn:GG:def}
		\hat \fX\;\vcentcolon=\;\tilde \fX+\perturbation\fI_{\stateDim} .
	\end{equation}
	satisfies the \LMI \eqref{eqn:LMI:P}.
\end{theorem}
Before proceeding with the proof of \Cref{teo4}, let us discuss the restriction on the type of \SLS imposed by assumption \eqref{eqn:ass:ndeAj}. By the problem setup, we have that $\fA_j$ is Hurwitz; hence, it immediately follows that if $\fA_j$ is normal, then \eqref{eqn:ass:ndeAj} holds since $\lambda_i(\fA_j+\fA_j^{\T})=2\real(\lambda_i(\fA_j))$. Also, it is known that every real Hurwitz matrix can be factorized as
\begin{equation}
	\label{eqn:dH}
	\fA_j=(\fJ_j-\fR_j)\fQ_j,
\end{equation}
with $\fJ_j\in\R^{\stateDim\times\stateDim}$ skew-symmetric, $\fR_j\in\R^{\stateDim\times\stateDim}$ symmetric positive semi-definite, and $\fQ_j\in\R^{\stateDim\times\stateDim}$ symmetric positive definite; cf.~\cite{GilS17}. Whenever $\fQ_j=\fI_{\stateDim}$ for all $j\in\switchingSet$, it is straightforward to check that \eqref{eqn:ass:ndeAj} is verified. Moreover, if the \SLS is quadratically stable, then there exists a symmetric positive-definite $\fQ$ such that one can choose $\fQ_j = \fQ$ for all $j\in\switchingSet$ in~\eqref{eqn:dH}, such that~\eqref{eqn:ass:ndeAj} is satisfied after a suitable coordinate transformation; see \cite{SchU18}.
\begin{proof}[Proof of~\Cref{teo4}]
	By the assumptions, $\fX$ satisfies the \LMIs 
	\begin{equation}\label{eqn:1:teo4}
	\fA_j \fX+\fX\fA_j^{\T}+\fB_j\fB_j^\T \prec \zeroVec \quad \text{for all } j\in\switchingSet,
  \end{equation}
  and we need to show that
  	\begin{equation}\label{eqn:2:teo4}
  \fw^{\T}\left(	\fA_j \hat \fX+\hat \fX\fA_j^{\T}+\fB_j\fB_j^\T \right)\fw< \zeroVec \quad \text{for all } j\in\switchingSet \text{ and for all }\fw\in\R^{\stateDim}\setminus \{\zeroVec\}.
  \end{equation}
  Using \eqref{eqn:1:teo4}, we observe $\fw^{\T}\fB_j\fB_j^\T\fw<-2\fw^{\T}\fX\fA^{\T}_j\fw$. Substituting into \eqref{eqn:2:teo4} yields
  \begin{equation}\label{eqn:3:teo4}
  2\fw^{\T}\hat \fX\fA_j^{\T}\fw+\fw^{\T}\fB_j\fB_j^\T \fw<2\fw^{\T}\hat \fE\fA_j^{\T}\fw
  \end{equation}
  with $\hat \fE\vcentcolon=\hat \fX-\fX$. Defining $\tilde{\fE} \vcentcolon= \tilde{\fX} - \fX$ and using~\eqref{eqn:GG:def} implies
  \begin{align*}
  	\begin{aligned}
  \fw^{\T}\hat \fE\fA_j^{\T}\fw= \fw^{\T}\tilde \fE\fA_j^{\T}\fw+\perturbation\fw^{\T}\fA_j^{\T}\fw\le\|\fw\|_2^2\left(\sigma_1({
  \fA_j})\|\tilde \fE\|_2-\sigma_\stateDim(\fA_j)\perturbation \right),
    	\end{aligned}
  \end{align*}
  where we used the fact that $\fw^{\T}\fA_j^{\T}\fw<0$ for all $\fw\in\R^{\stateDim}\setminus \{\zeroVec\}$ due to \eqref{eqn:ass:ndeAj}. The result thus follows immediately from~\eqref{eqn:ass:teo4}.
   \end{proof}
   
\begin{remark}
	The statement of~\Cref{teo4} can also be formulated in terms of the observability Gramian, in which case we need to assume that the \LMIs~\eqref{eqn:LMI:Q} are satisfied. The proof is analogous to the proof above.
\end{remark}

   Therefore, if we compute $\tilde \fX$ such that $ \|\fX-\tilde \fX\|_2\le\tol$, then we can choose 
   \begin{equation}\label{eqn:cond:LMI}
   	\perturbation=\hat \perturbation\vcentcolon=\tol \max_{j=1,\ldots,M}\kappa(\fA_j)
   \end{equation} 
   to enforce the \LMIs for the approximated solution. However, for large-scale systems arising from some \PDE discretization, $\kappa(\fA_j)$ may be large, resulting in $\perturbation$ being several orders of magnitude larger than $\tol$. In this way, the sum of the neglected singular values of $\hat \fX$ can't be less than $\perturbation$, dramatically affecting the sharpness of the error estimate \eqref{eq:err:est}. This consideration also led us to introduce the novel \MOR approach discussed in the next section.
\section{\MOR for \SLS using piecewise constant Gramians}\label{sec:PMOR}
The primary difficulty in \MOR for \SLS, when employing approximate \GLEs solutions, lies in simultaneously guaranteeing the stability of the resulting \ROM and the applicability of the error bound stated in~\Cref{thm:BTerrorBound}. Achieving this requires that the solutions of the \GLEs in \eqref{eqn:GLE} satisfy the \LMIs~\eqref{eqn:LMI}. For a specific class of \SLS, \Cref{teo4} provides a way to perturb an approximate \GLE solution so that the resulting perturbed solution satisfies the \LMIs, assuming the exact solution does as well. However, the condition in \eqref{eqn:cond:LMI} within \Cref{teo4} can reduce the sharpness of the error bound and applies only to \SLS that satisfy \eqref{eqn:ass:ndeAj}. Furthermore, it may not always be clear whether the exact solutions of the \GLEs satisfy the \LMIs in the first place. To overcome these limitations, we propose a novel \MOR strategy for \SLS that guarantees the stability of the \ROM and provides a generalized error bound without requiring the \LMIs to hold.

\subsection{Piecewise balancing reduction for \SLS}
\label{subsec:PBR}
We start with the observation that if in~\eqref{eqn:GLE:reach} we set $\fA=\fA_j$, then it is straightforward to see that the \GLE solution satisfies the \LMI
\begin{equation}\label{eqn:GLE:ineq}
	\fA_j\fX+\fX\fA_j^\T+\fB_j\fB_j \preceq 0.
\end{equation}
Therefore, each exact solution of the \GLE with $\fA=\fA_j$ always satisfies at least the \LMI~\eqref{eqn:GLE:ineq}. This inspires the following idea: for each $j\in\switchingSet$, solve the \GLE with $\fA=\fA_j$ (and correspondingly updated $\fN_i$), and, following the \MOR procedure described in \Cref{sec:MORexplained}, compute the projection matrices \eqref{eq:PGproj:mat}, which we denote as $\fV_j$ and $\fW_j$ to highlight that they originate from the \GLE where $\fA=\fA_j$. Then, we introduce the following reduced systems with state transition matrices
\begin{equation}\label{eqn:new:red:system}
\switchedSysRedJumps \quad \left\{\quad \begin{aligned}
		\dot{\stateRed}(t)  &= \reduceBis{\fA}_{\switch(t)} \stateRed(t) + \reduceBis{\fB}_{\switch(t)}\inp(t), &	\stateRed(t_0) &= \zeroVec, \\
		\stateRed(t^{+}_{k})&=\fW^{\T}_{\switch_k}\fV_{\switch_{k-1}}\stateRed(t^-_k),\\
		\reduce{\out}(t) &= \reduceBis{\fC}_{\switch(t)}\stateRed(t),\\
	\end{aligned}\right.
\end{equation}
where $t_k$ are the switching times for the given signal $\switch\in\calS$, $\switch_k\vcentcolon=\switch(t_k)$, and 
\begin{align}\label{eqn:new:PGproj:mat}
	\reduceBis{\fA}_j &\vcentcolon= \fW_j^\T \fA_j \fV_j \in \R^{\stateDimRed_j\times \stateDimRed_j}, &
	\reduceBis{\fB}_j &\vcentcolon= \fW_j^\T \fB_j\in \R^{\stateDimRed_j\times \inpDim}, &
	\reduceBis{\fC}_j &\vcentcolon= \fC_j \fV_j\in\R^{\outDim\times \stateDimRed_j}
\end{align}
for $j\in\switchingSet$. 
We call this procedure the \emph{piecewise balancing reduction} (\PBR) for \SLS. In the next sections, we discuss the stability of the reduced \SLS \eqref{eqn:new:red:system}. 

\subsection{Stability of the \PBR for \SLS}
We start with the following proposition.
\begin{proposition}\label{prop:stab}
Assume that the \GLEs~\eqref{eqn:GLE} are solved exactly with $\fA = \fA_j$ and correspondingly updated $\fN_i$ for $i,j\in\switchingSet$. Then, each subsystem of the switched system $\switchedSysRedJumps $ defined in \eqref{eqn:new:red:system} is stable.
\end{proposition}
\begin{proof}
	The stability arises from the observation that the solutions $\fcP_j$ and $\fcQ_j$ of \GLEs can be interpreted as generalized Gramians for the mode $j$ of the \SLS; see \cite[Sec.~4.7]{DulP00}. Therefore, each $\tilde \fA_j$ is, by construction, stable due to the stability properties of classical balance truncation.
\end{proof}
As discussed in \Cref{sec:MORexplained}, the fact that each subsystem mode is stable is not sufficient to conclude that the \SLS is stable for each $\switch\in \calS$. However, having stable subsystems represents an improvement over what the reduction procedure described in \Cref{sec:MORexplained} can achieve. Indeed, for that strategy, stable subsystems cannot be guaranteed if at least one of the \LMIs in \eqref{eqn:LMI} is not satisfied by one of the exact solutions of \eqref{eqn:GLE}. 

To further investigate the stability of \eqref{eqn:new:red:system}, one can attempt to solve the following set of \LMI 
\begin{equation*}
	\reduceBis{\fA}_j^{\T}	\reduceBis{\fP}+	\reduceBis{\fP}	\reduceBis{\fA}_j \prec \zeroMat\quad\text{for all }j\in\switchingSet,
\end{equation*}
to verify the quadratic stability property for the \ROM. Since the \SLS is now of size $\stateDimRed\ll\stateDim$, solving the \LMI is a much cheaper computational task. Moreover, if the solution of the \LMI is not successful, it is possible to exploit the fact that the asymptotic stability of the \SLS can still be guaranteed under a minimum switching time condition; see \cite[Sec.~3.2.1]{Lib03}. 

\subsection{Accuracy of the \PBR for \SLS}\label{sec:4.2}
To obtain an error bound for the \PBR, we will rely on standard techniques used to derive the classical balanced truncation error estimate. In particular, we start by investigating the effect of truncating only the smallest singular value; see the forthcoming \Cref{teo5}, and then iteratively apply this result to obtain the error bound, which we present in the forthcoming \Cref{teo6}. In contrast to the classical setting, however, we must take into consideration that the Gramians are time-dependent (piecewise constant in time) and the possibility that the \LMIs may not be exactly satisfied. Let us start by defining the following piecewise constant matrices
\begin{equation}\label{eqn:picewise:Gra}
	\tilde \fcP(t)\;\vcentcolon=\tilde\fcP_{\switch(t)},\quad\quad\tilde \fcQ(t)\;\vcentcolon=\tilde \fcQ_{\switch(t)},
\end{equation} 
where $\tilde \fcP_j$ and $\tilde \fcQ_j$ are the approximated solutions of \eqref{eqn:GLE:reach} and \eqref{eqn:GLE:observ}, respectively, for $\fA=\fA_j$ obtained via \Cref{alg:statIterGLE}. Then, for each $j\in\switchingSet$, consider the factorizations $\tilde \fcP_j=\tilde \fS_j\tilde \fS_j^{\T}$, $\tilde \fcQ_j=\tilde \fR_j\tilde \fR_j^{\T}$, and $\tilde \fS_j^{\T}\tilde \fR_j=\fU_j\fSigma_j\fV_j^{\T}$. The transformations to balance mode $j$ are given by 
\begin{equation*}
	\bar \fV_j\;\vcentcolon=\tilde \fS_j\fU_j\fSigma_j^{-1/2}\quad\text{and}\quad\bar\fW_j\;\vcentcolon=\;\tilde \fR_j\fV_j\fSigma_j^{-1/2} =\bar \fV^{-\T}_j.
\end{equation*}
Balancing each mode of system \eqref{eqn:sDAE} and performing the truncation, we obtain the \SLS~\eqref{eqn:new:red:system} with state transitions. For a switching signal $\switch\in\calS$, let us introduce the set of switching times up to time $t$, which we denote as $\calT_{\switch(t)}\vcentcolon=\{t_1,\ldots,t_K\}$, with $K\in\N$. Now, let us suppose that each subsystem $j$ of the \SLS in \eqref{eqn:sDAE} is balanced via $\tilde \fcP_j$ and $\tilde \fcQ_j$; thus, the generalized observability and controllability Gramians are of the form
\begin{equation}\label{eqn:bal:gra}
	\fSigma_j=\begin{bmatrix}
		\fSigma_{j,1}&\zeroVec\\
		\zeroVec&\sigma(j)\fI_{m(j)}
	\end{bmatrix} \quad\text{with}\quad\fSigma_{j,1}\in\R^{(\stateDim-m(j))\times(\stateDim-m(j))}\;\; \text{diagonal},
\end{equation}
and where $m(j)$ is the multiplicity of $\sigma(j)$, implying that $\sigma(j)$ is not an entry in $\fSigma_{j,1}$. For simplicity, let us assume $\tilde m=m(j)$ for all $j\in\switchingSet$ (similar results can be stated  considering the minimum value of the multiplicity among the $M$ system modes).  Because of \eqref{eqn:bal:gra}, it is convenient for us to rewrite the system matrices $\fA_j$, $\fB_j$, and $\fC_j$ for each $j\in\switchingSet$ in their balanced form
\begin{equation}\label{eqn:SLS:balanced}
	\bar{\fA}_{j}\vcentcolon=\begin{bmatrix}
		\reduceBis{\fA}_{j}& \fA_{j,12}\\
		\fA_{j,21}& \fA_{j,22}
	\end{bmatrix},\quad  \bar{\fB}_{j}\vcentcolon=\begin{bmatrix}
		\reduceBis{\fB}_{j}\\
		\fB_{j,2}
	\end{bmatrix},\quad  \bar{\fC}_{j}\vcentcolon=\begin{bmatrix}
		\reduceBis{\fC}_{j}&
		\fC_{j,2}
	\end{bmatrix},
\end{equation}
where $\reduceBis{\fA}_{j}, \reduceBis{\fB}_{j}$, and $\reduceBis{\fC}_{j}$ form the $j$th subsystem of  \eqref{eqn:new:red:system} with $\stateDimRed=\stateDim-\tilde{m}$. 

\begin{lemma}\label{teo5}
Consider the \ROM~\eqref{eqn:new:red:system} for $\stateDimRed=\stateDim-\tilde m$ and the \FOM in balanced form \eqref{eqn:SLS:balanced}. Take the decomposition $\stateBal(t) = [\stateBal_1(t)^\T, \stateBal_2(t)^\T]^\T$ where $\stateBal_1(t) \in \R^{\stateDimRed}$ and let us define the vectors
\begin{equation}\label{eqn:new:var}
	\bfX_c(t)\vcentcolon=\begin{bmatrix}
		\stateBal_1(t)+\stateRed(t)\\
		\stateBal_2(t)
	\end{bmatrix}\quad\text{and}\quad\bfX_o(t)\vcentcolon=\begin{bmatrix}
		\stateBal_1(t)-\stateRed(t)\\
		\stateBal_2(t)
	\end{bmatrix}.
\end{equation}
For given $\switch\in\calS$, with switching times $t_k$ ($k=1,\ldots,K$), we introduce the quantities
\begin{align}\label{eqn:swi:con:fac}
	\begin{aligned}
	\gleG_{o,\switch}(\stateDim,\stateDimRed)\;\vcentcolon=&\;{\bfX}_o(t)^{\T}\fSigma_{\switch(t)}{\bfX}_o(t) +\sum_{k=1}^{K}	\left({\bfX}_o(t^{-}_k)^{\T}\fSigma_{\switch(t^{-}_k)}{\bfX}_o(t^{-}_k)-{\bfX}_o(t^{+}_k)^{\T}\fSigma_{\switch(t^{+}_k)}{\bfX}_o(t^{+}_k)    \right),\\
	\gleG_{c,\switch}(\stateDim,\stateDimRed)\;\vcentcolon=&\;	\sigma(\switch(t_K^{+}))^2{\bfX}_c(t)^{\T}\fSigma^{-1}_{\switch(t)}{\bfX}_c(t) \\
	+&\;\sum_{k=1}^{K}	\left(\sigma(\switch(t^{-}_k))^2{\bfX}_c(t^{-}_k)^{\T}\fSigma^{-1}_{\switch(t^{-}_k)}{\bfX}_c(t^{-}_k)-\sigma(\switch(t_k^{+}))^2{\bfX}_c(t^{+}_k)^{\T}\fSigma^{-1}_{\switch(t^{+}_k)}{\bfX}_c(t^{+}_k)    \right),\\
	\fM_c(j)\;\vcentcolon=&\;\bar \fA_j\Sigma_j+\Sigma_j\bar \fA_j^{\T}+\bar \fB_j\bar\fB^{\T}_j,\quad	\fM_o(j)\;\vcentcolon=\;\bar \fA^{\T}_j\Sigma_j+\Sigma_j\bar \fA_j+\bar \fC_j^{\T}\bar\fC_j,\\
	\gleH_{c,\switch}(\stateDim,\stateDimRed)\;\vcentcolon=&\;\sum_{k=0}^{K}\lambda_{1}(	\fM_c(\switch(t_k^{+})))\int_{t_k}^{t_{k+1}}\|\bfX_c(s)\|_2^2\,\ds,\\
		\gleH_{o,\switch}(\stateDim,\stateDimRed)\;\vcentcolon=&\;\sum_{k=0}^{K}\lambda_{1}(	\fM_o(\switch(t_k^{+})))\int_{t_k}^{t_{k+1}}\|\bfX_o(s)\|_2^2\,\ds.
	\end{aligned}
\end{align}
 Then,  the output between the \FOM~\eqref{eqn:sDAE} and the \ROM~\eqref{eqn:new:red:system} of size $\stateDimRed=\stateDim-\tilde{m}$ satisfies the relation
\begin{align}
	\label{eq:err:est:2}
	\begin{aligned}
\left(\int_{0}^{t}\|\out(s)-\tilde{\out}(s)\|^2_{2}\,\ds\right)^{\frac{1}{2}} \;\leq &\;\Big(\gleH_{c,\switch}(\stateDim,\stateDimRed)+\gleH_{o,\switch}(\stateDim,\stateDimRed)-\left(\gleG_{o,\switch}(\stateDim,\stateDimRed)+\gleG_{c,\switch}(\stateDim,\stateDimRed)\right)\\
+&\;4 \tilde{\sigma}^2\int_{0}^{t}\|\inp(s)\|^2_{2}\,\ds\Big)^{\frac{1}{2}},
	\end{aligned}
\end{align}
with $\tilde \sigma\vcentcolon=\max_{j=1,\ldots,M}\sigma(j)$
\end{lemma}
\begin{proof}
	The proof is provided in \Cref{A-A1}. It follows the steps of \cite[Lem.~13]{PetWL13} and \cite{San04}, although these references do not consider piecewise Gramians and assume that the \LMI are satisfied.
\end{proof}

Before iteratively applying \Cref{teo5} to obtain the error bound, let us introduce and comment on some notation. Given $\stateDimRed\ll\stateDim$ and $k\in\N$ with $k>\stateDimRed$, we denote the distinct singular values of the balanced Gramian \eqref{eqn:bal:gra} and their associated multiplicities by $\sigma_k(j)$ and $m_k(j)$, for  $k= \stateDimRed+1, \ldots ,\tilde \stateDim$, respectively. In order to simplify the notation, we impose the assumption that $ m_{ k}(j_1)= m_k(j_2)$ holds for all $j_1,j_2\in\switchingSet$ and for all $k= \stateDimRed+1, \ldots ,\tilde \stateDim$. This allows us to drop the $j$ dependence for $m_k(j)$ and avoids the introduction of additional heavy notation in the statement of the next theorem. In general, we could state the next results without distinguishing between simple and non-simple singular values. However, as we will comment later, it is crucial to deal with the case where the smallest singular value of the Gramian \eqref{eqn:bal:gra} is not simple, as this allows for a sharper error bound that is also possibly efficiently computable. 

\begin{theorem}[Error bound of the \PBR for \SLS]\label{teo6}
	Given a switching signal $\switch\in\calS$, consider the \ROM~\eqref{eqn:new:red:system} of size $\stateDimRed$. Then, the output error between the \FOM~\eqref{eqn:sDAE} and the \ROM~\eqref{eqn:new:red:system} of size $\stateDimRed$ satisfies the relation 
	\begin{align}
		\label{eq:err:est:3}
		\begin{aligned}
	&\left(\int_{0}^{t}\|\out(s)-\tilde{\out}(s)\|^2_{2}\,\ds\right)^{\frac{1}{2}}
	\;\leq\;2\left(\int_{0}^{t}\|\inp(s)\|^2_{2}\,\ds\right)^{\frac{1}{2}}\sum_{k=\stateDimRed+1}^{\tilde \stateDim}\tilde{\sigma}_{s_{k}}\\
	+&
	\sum_{k=\stateDimRed+1}^{\tilde \stateDim}   \gleH_{c,o,\switch}(s_{k+1},s_{k})^{\frac{1}{2}}H(   \gleH_{c,o,\switch}(s_{k+1},s_{k}))\\
	-&
	\sum_{k=\stateDimRed+1}^{\tilde \stateDim}   \gleG_{c,o,\switch}(s_{k+1},s_{k})^{\frac{1}{2}}H(-   \gleG_{c,o,\switch}(s_{k+1},s_{k})),
			\end{aligned}
	\end{align}
	where $H(x)$ is the Heaviside function, $\tilde{\sigma}_k\vcentcolon=\max_{j\in\switchingSet}\sigma_{k}(j) $, $s_k\vcentcolon=\sum_{k=\stateDimRed}^{\tilde n}m_k$, $m_\stateDimRed\vcentcolon=\stateDimRed$, and
	\begin{align*}
	     \begin{aligned}
	     \gleH_{c,o,\switch}(s_{k+1},s_{k})	\;&\vcentcolon=\;\gleH_{o,\switch}(s_{k+1},s_{k})+\gleH_{c,\switch}(s_{k+1},s_{k}),\\
	       \gleG_{c,o,\switch}(s_{k+1},s_{k})	\;&\vcentcolon=\;\gleG_{o,\switch}(s_{k+1},s_{k})+\gleG_{c,\switch}(s_{k+1},s_{k}).
	     \end{aligned}
	\end{align*}
\end{theorem}
\begin{proof}
The proof follows directly by applying the sub-additivity of the square root function to \Cref{teo5} and standard balance truncation error arguments; see, for instance, \cite[Sec.~7]{PetWL13}.
\end{proof}
The error bound presented in \Cref{teo6} comprises three distinct contributions: the first corresponds to the sum of the truncated singular values, as in \Cref{thm:BTerrorBound}; the second stems from the use of piecewise constant Gramians in the derivation of \eqref{eqn:new:red:system}; and the third reflects the fact that the numerically computed solutions of the \GLE may not exactly satisfy the \LMI conditions due to approximation errors. In this regard, \Cref{teo6} extends \Cref{thm:BTerrorBound} by incorporating the effects of piecewise constant Gramians and does not require that the \LMI conditions be fulfilled exactly. However, this generalization comes at the cost of requiring the reduced state trajectory to be computed in order to evaluate the error bound. This a-posteriori evaluation step is not required for the bound in \Cref{thm:BTerrorBound}.

\begin{remark}\label{rmk:GLEvsLE}
	The proof of \Cref{teo5} requires \eqref{eqn:bal:gra} to be invertible. To maintain computational efficiency, we apply truncation when computing the approximate \GLE solution 
	$\tilde{\fX}$, i.e., we discard all singular values smaller than the prescribed accuracy threshold 
	$\tol$. As a result, \eqref{eqn:bal:gra} becomes singular once the reduced dimension is sufficiently large. To address this, we replace the zero singular values with 
	$\tol$, thereby ensuring that \eqref{eqn:bal:gra} remains invertible. Moreover, in terms of the error estimates, since we are introducing non-simple singular values, we only need to take them into account once in the error estimate \eqref{teo6}. In practice, this is done by adding $\tol$ to each of the diagonal entries of $\Sigma_j$ for $j=1,\ldots,M$.
\end{remark}

\begin{remark}
   The rationale for employing the solution of the \GLE for the \PBR instead of the simpler solutions of standard Lyapunov equations is that the latter do not provide sufficient information about the reachable and observable sets of the full system; that is, they do not satisfy \eqref{eqn:reach:observ:gramians}. As a direct consequence, the jump terms appearing in the corresponding error bound can, in general, be expected to be larger when \eqref{eqn:new:red:system} is constructed via simple Lyapunov equations; see \Cref{fig6}.
\end{remark}

\subsection{Tightening the error bound}\label{subsec:4.3}
In the following, we propose a methodology to tighten the error bound \eqref{eq:err:est:3} enforcing the non-positivity of the eigenvalues of $\fM_c(j)$ (and $\fM_o(j)$) by suitably perturbing $\tilde{\fX}_j$. Plugging $\tilde \fX_j$ into the \GLE for $\fA=\fA_j$ with correspondingly upgraded $\fN_j$, we can define the $j$-th residual
\begin{equation*}
 \tilde \fR_j\;\vcentcolon=\;	\fA_j \tilde \fX_j+\tilde \fX_j\fA_j^\T+\sum_{i=1}^{M}\left(\fN_i \tilde \fX_j\fN_i^\T+\fB_i\fB_i^\T \right).
\end{equation*}
The matrix $ \tilde \fR_j$ is symmetric and, in general, indefinite. Let us write its eigenvalue decomposition as
\begin{equation}\label{eqn:EIG:DEC:RES}
 \tilde \fR_j\;=\;\begin{bmatrix}
 	\fV^{-}_j& 	\fV^{+}_j
 \end{bmatrix}\begin{bmatrix}
 \fSigma^{-}_j& \zeroMat\\
 \zeroMat &  \fSigma^{+}_j
 \end{bmatrix}\begin{bmatrix}
 \fV^{-}_j& 	\fV^{+}_j
 \end{bmatrix}^{\T},
\end{equation} 
where $ \fSigma^{-}_j$, $ \fSigma^{+}_j$  are the diagonal matrices containing the non-positive and the positive eigenvalues, respectively, while the columns of the matrices  $\fV^{-}_j$, $\fV^{+}_j$ are the associated orthonormal eigenvectors. We now define $\fDelta_j$ as the solution to the Lyapunov equation
\begin{equation}\label{eqn:pert:lyap}
	\fA_j\fDelta_j+	\fDelta_j\fA_j^{\T}+	\fV^{+}_j  \fSigma^{+}_j(\fV^{+}_j)^{\T}\;=\;\zeroVec .
\end{equation}
\begin{proposition}\label{prop5}
	The matrix $\hat \fX_j\vcentcolon=\tilde{\fX}_j+\fDelta_j$ with $\fDelta_j$ being the unique solution of~\eqref{eqn:pert:lyap} satisfies 
	\begin{equation*}
			\fA_j\hat \fX_j+\hat \fX_j\fA_j^{\T}+	\fB_j\fB_j^{\T} \;\preceq\;\zeroVec.
	\end{equation*}
\end{proposition}
\begin{proof}
	We have
 \begin{align*}
 	\begin{aligned}
 			\fA_j\hat \fX_j+\hat \fX_j\fA_j^{\T}+	\fB_j\fB_j^{\T} \;\preceq&\; 	\fA_j\hat \fX_j+\hat \fX_j\fA_j^{\T}+\sum_{i=1}^{M}\left(\fN_i \tilde \fX_j\fN_i^\T+\fB_i\fB_i^\T \right)\\
 			=&\;\tilde{\fR}_j+	\fA_j\fDelta_j+	\fDelta_j\fA_j^{\T}\;=\;	\fV^{-}_j  \fSigma^{-}_j(\fV^{-}_j)^{\T}\;\preceq\;\zeroMat, 
 			\end{aligned}
 \end{align*}
 where we used \eqref{eqn:pert:lyap} and \eqref{eqn:EIG:DEC:RES}.
\end{proof}

\subsection{Numerical aspects concerning \texorpdfstring{\Cref{sec:4.2,subsec:4.3}}{Sections~4.3 and~4.4}}\label{subsec:NumAspectsSec4}

We start detailing how the numerics involved in computing and storing $\hat \fX_j$ can be performed efficiently. 

\begin{enumerate}
\item For any $j\in\switchingSet$, we need to compute the positive eigenvalues and associated eigenvectors of $\tilde \fR_j$. Their computation can be performed without explicitly forming~$\tilde \fR_j$, exploiting the low-rank factor $\tilde{\fZ_j}$ of $\tilde{\fX}_j$ and employing only matrix-vector products; see \cite{LehSY98}.
\item For each $j\in\switchingSet$, our method requires solving the Lyapunov equation \eqref{eqn:pert:lyap} in addition to using \Cref{alg:statIterGLE} to compute $\tilde{\fX}_j$. Solving Lyapunov equations still introduces numerical errors; however, being a simpler matrix equation than the \GLE, we can solve it with higher accuracy without compromising the overall computational cost.
\item The evaluation of $\fDelta_j$ is more efficient as fewer positive eigenvalues are present in $\tilde \fR_j$ for all $j\in\switchingSet$. Their number may not be small in general; however, even computing a few of the largest ones is typically enough to tighten the error bound and improve the accuracy of the \ROM.
\item We do not explicitly compute $\fDelta_j$, but rather compute an approximation $\tilde \fDelta_j$ by storing its low-rank factor, say $\tilde \fG_j$, such that $\tilde \fDelta_j = \tilde \fG_j\tilde \fG_j^{\T}$. In particular, \eqref{eqn:pert:lyap} must be solved with significantly higher accuracy than the corresponding \GLE. To form $\hat \fX_j$, we append the low rank factors of $\tilde \fDelta_j$ to those of $\tilde \fX_j$, i.e.,
\begin{equation}\label{eqn:mod:Gra}
\hat \fX_j\;=\;\begin{bmatrix}
\tilde \fZ_j& \tilde \fG_j
\end{bmatrix}\begin{bmatrix}
\tilde \fZ_j& \tilde \fG_j
\end{bmatrix}^{\T}.
\end{equation}
\end{enumerate}
We conclude by commenting on the error bound of \Cref{teo6}. Evaluating \Cref{teo6} requires computing all reduced state trajectories from dimension $\stateDimRed$ to $\stateDim$, which is computationally inefficient. However, since the reduced system is in a balanced form, the entries of the state decay rapidly. Therefore, it is reasonable to approximate the reduced state at higher dimensions using the trajectory computed at a smaller dimension—provided that the last reduced state components are significantly smaller than the target accuracy.
\section{Numerical Experiments}
\label{sec:examples}
In this section, we present numerical experiments to support both the theoretical results and the proposed methodology. We begin by validating the theory developed in \Cref{sec:numericalDetails} using a synthetic example from \cite[Sec.~V]{PonGB20}. Subsequently, we illustrate the theoretical framework of \Cref{sec:PMOR} through a test case based on the parametric Black-Scholes model, where switching behavior is induced by treating the parameter as piecewise constant in time. All calculations were performed with MATLAB 2024b on a laptop with an Apple M2 Pro processor and 16GB of RAM. 

\vspace{0.2cm}
\noindent\fbox{%
	\parbox{0.98\textwidth}{%
		The code and data used to generate the subsequent results are accessible via
		\begin{center}
			\url{https://zenodo.org/records/19923899}
		\end{center}
		under MIT Common License.
	}%
}\\[.2em]
\subsection{Synthetic example} 
\label{subsec:MSD}
Let us consider a 2-mode switched linear system of dimension~$\stateDim$, with state matrices given by $\fA_1, \fA_2 \in \R^{\stateDim\times\stateDim}$ such that
\begin{align}\label{eqn:synt:ex}
	\begin{aligned}
    \fA_1(i,j)\vcentcolon=\begin{cases} -1, & \mbox{if }i=j \\ 
    	\frac{1}{2}, & \mbox{if }i-1=j\\
    	0,&\mbox{otherwise}
    \end{cases},\quad\quad \fA_2(i,j)\vcentcolon=\begin{cases} -2, & \mbox{if }i=j \\ 
    \frac{4}{5}, & \mbox{if }i-1=j\\
    -\frac{1}{5}, & \mbox{if }i+1=j\\
    0,&\mbox{elsewhere}
    \end{cases}
	\end{aligned}
\end{align}
and $\fB_1^\T=\fC_1=[0,\ldots,0,1]$, $\fB_2^\T=\fC_2=[1,0,\ldots,0]$. We set $\fA =\fA_1$, $\fN_1=\zeroVec$, and $\fN_2 =\fA_2 - \fA_1$ to obtain the \GLE~setting in~\eqref{eqn:generic:GLE}. Note that both $\fA_1$ and $\fA_2$ are strictly dissipative matrices.

\begin{figure}[t]
	\centering
	\begin{subfigure}[t]{.48\linewidth}
%
\begin{tikzpicture}

\begin{axis}[%
	width=\imageWidth,
	height=\imageHeight,
	scale only axis,
	scaled ticks=false,
	grid=both,
	grid style={line width=.1pt, draw=gray!10},
	major grid style={line width=.2pt,draw=gray!50},
	axis lines*=left,
	axis line style={line width=\lineWidth},
	xmode=log,
xmin=1e-12,
xmax=1e0,
xlabel style={font=\color{white!15!black}},
xlabel={$\tol$},,
ymode=log,
ymin=1e-12,
ymax=1,
yminorticks=true,
ylabel style={font=\color{white!15!black}},
ylabel={absolute error},
	axis background/.style={fill=white},
	legend style={%
		legend cell align=left, 
		align=left, 
		font=\tiny,
		draw=white!15!black,
		at={(0.30,0.65)},
		anchor=south,},
]

\addplot [color = mycolor1, line width=\lineWidth,mark=o, mark options={solid, mycolor1}] table [x=x1, y=y2, col sep=comma]{Img/DataCSV/Example_1_Plot_1.csv};
    \addlegendentry{$ \|\fcP-\tilde{\fcP}\|_2$}
    
    \addplot [color = mycolor2,dashed, line width=\lineWidth,mark=square, mark options={solid, mycolor2}] table [x=x1, y=y3, col sep=comma]{Img/DataCSV/Example_1_PLot_1.csv};
   \addlegendentry{$ \|\fcQ-\tilde{\fcQ}\|_2$}

        \addplot [color = mycolor3,dash dot, line width=\lineWidth, mark=triangle, mark options={solid, mycolor3}] table [x=x1, y=y1, col sep=comma]{Img/DataCSV/Example_1_PLot_1.csv};
   \addlegendentry{$\tol$}

\end{axis}
\end{tikzpicture}%
		\caption{Absolute error vs.~exit tolerance $\tol$ for the computing the Gramians with~\Cref{alg:statIterGLE} for $\stateDim=200$.}
		\label{fig1:left}
	\end{subfigure}
	\hfill
	\begin{subfigure}[t]{.48\linewidth}
%
\begin{tikzpicture}

\begin{axis}[%
width=\imageWidth,
height=\imageHeight,
at={(1.011in,0.642in)},
scale only axis,
grid=both,
grid style={line width=.1pt, draw=gray!10},
major grid style={line width=.2pt,draw=gray!50},
axis lines*=left,
axis line style={line width=\lineWidth},
mark size=2.2pt,
xmode=log,
xmin=50,
xmax=3e6,
xminorticks=true,
xlabel style={font=\color{white!15!black}},
xlabel={$\stateDim$},
ymode=log,
ymin=0.05,
ymax=1200,
yminorticks=true,
ylabel style={font=\color{white!15!black}},
ylabel={computational time (s)},
axis background/.style={fill=white},
legend style={%
	legend cell align=left, 
	align=left, 
	font=\tiny,
	draw=white!15!black,
	at={(0.30,0.65)},
	anchor=south,},
]
\addplot [color=mycolor1, line width=\lineWidth, mark=o, mark options={solid, mycolor1}]
  table [x=x1, y=y1, col sep=comma]{Img/DataCSV/Example_1_Plot_2.csv};
\addlegendentry{$\fcP$}

\addplot [color=mycolor2, line width=\lineWidth, mark=asterisk, mark options={solid, mycolor2}]
table [x=x1, y=y2, col sep=comma]{Img/DataCSV/Example_1_Plot_2.csv};
\addlegendentry{$\fcQ$}

\addplot [color=black, dashdotted, line width=\lineWidth]
table [x=x1, y=y3, col sep=comma]{Img/DataCSV/Example_1_Plot_2.csv};
\addlegendentry{$\mathcal{O}(\stateDim)$}

\end{axis}
\end{tikzpicture}%
		\caption{Computational time to solve the \GLEs with respect to the size of the full problem $\stateDim$.}
		\label{fig1:right}
	\end{subfigure}
	\caption{Synthetic example \eqref{eqn:synt:ex}.}
	\label{fig1}%
\end{figure}%
In our first experiment, we set $\stateDim = 200$, and we compare the approximate Gramians computed via \Cref{alg:statIterGLE} against the exact ones obtained by solving the full $\stateDim^2 \times \stateDim^2$ linear system associated with the operators in \eqref{eqn:GLE:oper}, for various values of $\tol$. For this problem, we have $\|\calL^{-1}\Pi\| = 0.8098$, indicating that no rescaling is necessary. The associated constant from~\eqref{eqn:GLEerrorBoundConstant} is given by $\gamma = 4.2576$. \Cref{fig1:left} illustrates the effectiveness of the error bound presented in \Cref{teo1} for computing the Gramians via \Cref{alg:statIterGLE}. Notably, the ratio between the prescribed tolerance $\tol$ and the observed error for both Gramians remains below 10, highlighting the sharpness of the error estimate used to define the stopping criteria in \Cref{alg:statIterGLE}. In \Cref{fig1:right}, we report the computational time (in seconds) required to execute \Cref{alg:statIterGLE} as a function of the problem dimension $\stateDim$. Thanks to the use of sparse matrices and the fact that \eqref{eqn:synt:ex} are dissipative matrices, the runtime exhibits asymptotic linear growth, i.e.,~$\mathcal{O}(n)$. Overall, \Cref{fig1} demonstrates that our implementation of \Cref{alg:statIterGLE} efficiently delivers accurate approximations for any user-specified tolerance $\tol$.

We observe that the exact Gramians of this example satisfy the set of \LMIs~\eqref{eqn:LMI} up to machine precision, indicating that the problem is well suited to ensure a quadratically stable \ROM. Moreover, both $\fA_1$ and $\fA_2$ are dissipative matrices, so the hypothesis of \Cref{teo4} holds. The spectral condition numbers are $\kappa(\fA_1)=3$ and $\kappa(\fA_2)=2.4136$. Let us define the operators $\gleM_i:\R^{\stateDim\times\stateDim}\rightarrow\R^{\stateDim\times\stateDim}$
\begin{equation}\label{eqn:def:M}
\gleM_i(\fX)\;\vcentcolon=\;\fA_i\fX+\fX\fA_i^\T+\fB_i\fB_i^{\T},\quad i=1,2.
\end{equation}
\begin{table}[t]
	\centering
	\begin{tabular}{ll@{\hspace{3em}}ll@{\hspace{3em}}ll}
		\toprule
		$\tol$ & $\hat \perturbation$ & $\lambda_1(\gleM_1(\tilde\calP))$ & $\lambda_1(\gleM_1(\hat \calP))$ & $\lambda_1(\gleM_2(\tilde\calP))$ & $\lambda_1(\gleM_2(\hat \calP))$  \\
		\midrule
		$10^{-2}$& 3.0$\;\cdot\;10^{-2}$ & 1.0$\;\cdot\; 10^{-3}$ & -3.0 $\cdot\;10^{-2}$ & 1.2$\;\cdot\; 10^{-3}$ & -8.4 $\cdot\;10^{-2}$ \\
		$10^{-4}$& 3.0$\;\cdot\;10^{-4}$ & 1.6$\;\cdot\; 10^{-5}$ & -3.0 $\cdot\;10^{-4}$ & 7.6$\;\cdot\; 10^{-6}$ & -8.4 $\cdot\;10^{-4}$ \\
		$10^{-6}$& 3.0$\;\cdot\;10^{-6}$ & 3.8$\;\cdot\; 10^{-8}$ & -3.0 $\cdot\;10^{-6}$ & 4.8$\;\cdot\; 10^{-8}$ & -8.4 $\cdot\;10^{-6}$ \\
		$10^{-8}$& 3.0$\;\cdot\;10^{-8}$ & 1.3$\;\cdot\; 10^{-9}$ & -3.0 $\cdot\;10^{-8}$ & 1.7$\;\cdot\; 10^{-9}$ & -8.4 $\cdot\;10^{-8}$ \\
		$10^{-10}$& 3.0$\;\cdot\;10^{-10}$ & 8.6$\;\cdot\; 10^{-12}$ & -3.0 $\cdot\;10^{-10}$ & 1.1$\;\cdot\; 10^{-11}$ & -8.4 $\cdot\;10^{-10}$ \\
		\bottomrule
	\end{tabular}
	\caption{Synthetic example \eqref{eqn:synt:ex} for $\stateDim=200$. The table reports the values of $\hat \perturbation$ (as defined in \eqref{eqn:cond:LMI}) and the satisfaction of the \LMIs in \eqref{eqn:LMI:P} for both the approximated Gramian $\tilde \fcP$ and the generalized Gramian $\hat \fcP$, across different values of the approximation tolerance $\tol$ used in solving the \GLE \eqref{eqn:GLE:reach}. 
	}
	\label{tab:1}
\end{table}

In \Cref{tab:1}, we report, for several values of the prescribed accuracy $\tol$, the computed $\hat\perturbation$ (as defined in \eqref{eqn:cond:LMI}) and the largest eigenvalue of the matrix in \eqref{eqn:def:M}, with $i=1,2$. These quantities are shown for both the approximated reachability Gramian $\tilde \fcP$ and its perturbed counterpart $\hat \fcP$, which is adjusted to satisfy the \LMIs, as described in \Cref{subsection:3.2.1}. While $\tilde \fcP$ consistently fails to satisfy the \LMIs, the perturbed version $\hat \fcP$ always fulfills them, thereby numerically validating \Cref{teo4}.

Let us denote the cumulative sum of the neglected singular values as
\begin{equation}\label{eqn:neg:SV}
\tau(\stateDimRed,\Upsilon)\;\vcentcolon=\;2\sum_{k=\stateDimRed+1}^{\stateDim} \sigma_k(\Upsilon),
\end{equation}
and
\begin{align}\label{eqn:sys:err}
	\eta(t)\;\vcentcolon=\;\left(\int_{0}^{t}\|\inp(s)\|_2^{2}\,\ds\right)^{\frac{1}{2}},\quad
	\varepsilon(\stateDimRed,\Upsilon)\;\vcentcolon=\;{\left(\int_{0}^{t}\|\out(s)-\tilde \out(s)\|^2_2\,\ds\right)^{\frac{1}{2}}}\eta(t)^{-1},
\end{align}
for $\Upsilon\in\{\tilde \fS^{\T} \tilde \fR,\,\hat \fS^{\T} \hat \fR \}$ where the output $\tilde \out$ is the one of the reduced system \eqref{eqn:sDAE:ROM} of size~$\stateDimRed$. For $\stateDim=2000$, we approximate the Gramians with $\tol=10^{-10}$ and display the results in \Cref{fig2}. The left plot shows the decay of the sum of the neglected singular values derived from the factors of the Gramian matrices. The green line corresponds to the exact Gramians, while the blue line represents the Gramians obtained by approximating the solution of the \GLE. As expected, the green curve lies above the blue one since $\fcP \succeq \tilde\fcP$ and $\fcQ \succeq \tilde\fcQ$. This property follows directly from an application of \cite[Thm.~3.3.16]{HorJ91}. The right plot compares the actual output error (see \eqref{eqn:sys:err}) with the error bound established in \Cref{thm:BTerrorBound}. When using the approximated Gramians $\tilde\fcP$ and $\tilde\fcQ$, the \LMI condition is not satisfied. As a result, the error bound becomes unreliable, especially as we approach the smallest neglected singular values. In contrast, by using the perturbed Gramians $\hat\fcP$ and $\hat\fcQ$, which are explicitly constructed to satisfy the \LMI, we obtain a \ROM for which the error estimate \eqref{eq:err:est} is rigorously verified.

\begin{figure}[t]	
	\centering{
		\begin{subfigure}[t]{0.48\textwidth}
%
\begin{tikzpicture}

\begin{axis}[%
width=\imageWidth,
height=\imageHeight,
at={(1.011in,0.642in)},
scale only axis,
grid=both,
grid style={line width=.1pt, draw=gray!10},
major grid style={line width=.2pt,draw=gray!50},
axis lines*=left,
axis line style={line width=\lineWidth},
mark size=2.2pt,
xmin=0,
xmax=31,
xlabel style={font=\color{white!15!black}},
xlabel={$\stateDimRed$},
ymode=log,
ymin=1e-16,
ymax=1e2,
yminorticks=true,
axis background/.style={fill=white},
legend style={%
	legend cell align=left, 
	align=left, 
	font=\tiny,
	draw=white!15!black,
	at={(0.99,0.99)},
	anchor=north east,},
]
\addplot [color=mycolor1, line width=\lineWidth, mark=o, mark options={solid, mycolor1}]
 table [x=x1, y=y1, col sep=comma]{Img/DataCSV/Example_1_Plot_3_T1.csv};
\addlegendentry{$\tau(\stateDimRed, \tilde 	\fS^\T\tilde \fR)$}

\addplot [color=mycolor2, line width=\lineWidth, mark=asterisk, mark options={solid, mycolor2}]
table [x=x2, y=y2, col sep=comma]{Img/DataCSV/Example_1_Plot_3_T2.csv};
\addlegendentry{$\tau(\stateDimRed,\hat 	\fS^\T\hat\fR)$ }
\addplot [color=mycolor4, line width=\lineWidth, mark=+, mark options={solid, mycolor4}]
table [x=x3, y=y3, col sep=comma]{Img/DataCSV/Example_1_Plot_3_T3.csv};
	\addlegendentry{$\tau(	\stateDimRed,\fS^\T\fR)$}
\end{axis}
\end{tikzpicture}%
			\subcaption{Decay of the sum of the neglected singular values, see \eqref{eqn:neg:SV}.}
			\label{fig2:a}
		\end{subfigure}
		\hfill
		\begin{subfigure}[t]{0.48\textwidth}
%
\begin{tikzpicture}

\begin{axis}[%
width=\imageWidth,
height=\imageHeight,
at={(1.011in,0.642in)},
scale only axis,
grid=both,
grid style={line width=.1pt, draw=gray!10},
major grid style={line width=.2pt,draw=gray!50},
axis lines*=left,
axis line style={line width=\lineWidth},
mark size=2.2pt,
xmin=0,
xmax=29,
xlabel style={font=\color{white!15!black}},
xlabel={$\stateDimRed$},
ymode=log,
ymin=0.5e-16,
ymax=1e3,
yminorticks=true,
axis background/.style={fill=white},
legend style={%
	legend cell align=left, 
	align=left, 
	font=\tiny,
	draw=white!15!black,
	at={(0.99,0.99)},
	anchor=north east,},
]
\addplot [color=mycolor1, line width=\lineWidth, mark=o, mark options={solid, mycolor1}]
 table [x=x1, y=y1, col sep=comma]{Img/DataCSV/Example_1_Plot_4.csv};
\addlegendentry{$\tau(\stateDimRed, \tilde \fS^{\T}\tilde \fR)$}

\addplot [color=mycolor2, line width=\lineWidth, mark=asterisk, mark options={solid, mycolor2}]
table [x=x1, y=y2, col sep=comma]{Img/DataCSV/Example_1_Plot_4.csv};
\addlegendentry{$\tau( \stateDimRed,\hat \fS^{\T} \hat \fR)$ }

\addplot [color=mycolor4, line width=\lineWidth, mark=+, mark options={solid, mycolor4}]
table [x=x1, y=y3, col sep=comma]{Img/DataCSV/Example_1_Plot_4.csv};
	\addlegendentry{$\varepsilon(\stateDimRed,\tilde \fS^{\T}\tilde \fR)$}

	\addplot [color=mycolor3, line width=\lineWidth, mark=square, mark options={solid, mycolor3}]
	table [x=x1, y=y4, col sep=comma]{Img/DataCSV/Example_1_Plot_4.csv};
	\addlegendentry{$\varepsilon(\stateDimRed, \hat \fS^{\T}\hat \fR)$}
	
\end{axis}
\end{tikzpicture}%
			\subcaption{\MOR error \eqref{eqn:sys:err} compared with the a priori error estimate of \Cref{thm:BTerrorBound}. Input $\inp=\sin(t^2+t)$.}
			\label{fig2:b}
		\end{subfigure}
		}
	\caption{Synthetic example \eqref{eqn:synt:ex} for $\stateDim=2000$. \MOR for \SLS through \GLEs solution.}
	\label{fig2}
\end{figure}

\subsection{Parametric Black-Scholes problem} 
\label{subsec:BlackSholes}
The well-known (deterministic) Black-Scholes equation (see, e.g., \cite{BlaS73}) has the following form:
\begin{equation}\label{eq:BS}
	\frac{\partial \psi}{\partial t}=\frac{1}{2}\nu^2s^2\frac{\partial^2\psi}{\partial s^2}+\ir s\frac{\partial \psi}{\partial s}-\ir \psi,\qquad s>L,\;\;0<t\le T,
\end{equation}
for $L$, $T$ given, where the unknown function $\psi(s,t)$ represents the fair price of the option when the corresponding asset price at time $T-t$ is $s$, and $T$ is the maturity time of the option. Moreover, $\ir\ge0$ and $\nu>0$ are given constants (representing the interest rate and the volatility, respectively). In practice, one considers a bounded spatial domain, setting $L<s<S$ for a sufficiently large $S$. We study \eqref{eq:BS} under the following conditions, which are a derivation from the European option call (see \cite{HouW07}) with initial data set to zero, i.e.,
\begin{equation}\label{ibcBS}
  \psi(L,t)=0, \quad \psi(S,t)=S-\ee^{-\ir t}K,\quad0\le t \le T.
\end{equation}
Following the same strategy adopted in \cite{HouW11}, we discretize the space on a uniform grid of $\stateDim =1000$ points in the interval $[0, 200]$ with $K=80$, using the classical centered finite difference scheme. We thus have the following discrete matrices
\begin{align}
	\label{eqn:ref:mat:BS}
	\fA(\nu,\ir) &= \frac{\nu^2}{2}\fD+\ir \fG, &
	\fB(\nu,\ir) &= \begin{bmatrix}
	S\nu^2\fB^\T_1+\ir\fB^\T_2\\
	-K\nu^2\fB^\T_1+\ir \fB^\T_2 
	\end{bmatrix}^\T, &
	\fC &= \begin{bmatrix}
		\frac{1}{\stateDim},\ldots,\frac{1}{\stateDim} \\
		0,\ldots,1
	\end{bmatrix},
\end{align}
where the output matrix $\fC$ accounts for the average fair price of the option and its value at the boundary $s=S$, $\fB_1$ and $\fB_2$ arise from the boundary conditions \eqref{ibcBS}. Furthermore, the matrices $\fD$ and $\fG$ correspond to the discretization of the diffusion and convection--reaction terms, respectively. To generate a switches system, we consider $M=4$ parameter samples $\fmu_i\vcentcolon=(\nu_i,\;\ir_i)$ given by
\begin{equation*}
	\fmu_1=[0.25,\;0.001],\quad	\fmu_2=[0.05,\;0.02],\quad	\fmu_3=[0.25,\;0.02],\quad	\fmu_4=[0.05,\;0.001].
\end{equation*}
We apply the \PBR method described in \Cref{sec:PMOR} to the resulting \SLS, setting $\tol=10^{-8}$ and normalizing the input matrices $\fB(\nu,\ir)$. 

First, we want to verify the effectiveness of the procedure described in \Cref{subsec:4.3} to tighten the \LMIs . With this aim, let us define the matrix operator $\gleN_i:\R^{\stateDim\times\stateDim}\rightarrow\R^{\stateDim\times\stateDim}$
\begin{equation}\label{eqn:def:N}
	\gleN_j(\fX)\;\vcentcolon=\;\fA_j^{\T}\fX+\fX\fA_j+\fC^{\T}\fC^{\T},\quad j=1,\ldots,4.
\end{equation}
\Cref{tab:2} displays the largest eigenvalue of $\gleN_j$ for each system mode. Using the modified Gramians \eqref{eqn:mod:Gra} successfully reduces the magnitude of the largest eigenvalues by several orders. 
\begin{table}[t]
	\centering
	\begin{tabular}{l@{\hspace{2em}}llll}
		\toprule
		  & $j=1$ &  $j=2$  &  $j=3$  &  $j=4$  \\ \midrule
		        \rule{0pt}{12pt}  $\lambda_{1}(\gleN_j(\tilde \fcQ_j))$  &$1.1\cdot10^{-5}$ & $2.9\cdot10^{-7}$ & $1.5\cdot10^{-5}$& $3.6\cdot10^{-7}$ \\ 
		        \rule{0pt}{12pt}  $\lambda_{1}(\gleN_j(\hat \fcQ_j))$   &$1.6\cdot10^{-11}$ &$6.8\cdot10^{-13}$  & $2.3\cdot10^{-11}$ & $5.6\cdot10^{-13}$ \\\bottomrule
	\end{tabular}
		\caption{Parametric Black-Scholes problem \eqref{eqn:ref:mat:BS}. The table shows the largest eigenvalue of~\eqref{eqn:def:N} for each system mode.}
	\label{tab:2}
\end{table}
Next, we compare the \FOM with the \ROM considering the input functions
\begin{equation}\label{eqn:inp:BS}
	\inp_1(t)\;\vcentcolon=\;\begin{bmatrix}
	1&\ee^{-0.02 t}
	\end{bmatrix}^{\T}\quad\text{and}\quad 	\inp_2(t)\;\vcentcolon=\;\begin{bmatrix}
\sin\left(2\pi\ee^\frac{t}{2}\right)&\sin\left(2\pi\ee^\frac{t}{2}\right)
	\end{bmatrix}^{\T}.
\end{equation}
Results for the two input functions are displayed in \Cref{fig4} for a switching signal with $10$ randomly chosen switching times. For both input functions, the $\ROM$ output with $\stateDimRed=30$ qualitatively overlaps with the \FOM. With the aim of showing the different error bound terms in \Cref{teo6}, we introduce the following functions:
\begin{align}\label{eqn:tau:tilde}
 \tau(\stateDimRed, \Upsilon)\;\vcentcolon=&\;2\sum_{k=\stateDimRed+1}^{\tilde \stateDim} \tilde \sigma_k,\\
\chi(\stateDimRed, \Upsilon)\;\vcentcolon=&\;	\eta(t)^{-1}\sum_{k=\stateDimRed+1}^{\tilde \stateDim}\gleH_{c,o,\switch}(s_{k+1},s_{k})^{\frac{1}{2}}H(\gleH_{c,o,\switch}(s_{k+1},s_{k})),\label{eqn:chi}\\
\iota(\stateDimRed, \Upsilon)\;\vcentcolon=&\;-
	\eta(t)^{-1}\sum_{k=\stateDimRed+1}^{\tilde \stateDim}\gleG_{c,o,\switch}(s_{k+1},s_{k})^{\frac{1}{2}}H(-\gleG_{c,o,\switch}(s_{k+1},s_{k})),\label{eqn:iota}\\
	 \varphi(\stateDimRed ,\Upsilon)\;\vcentcolon=&\; \tau(\stateDimRed; \Upsilon)+ \chi(\stateDimRed; \Upsilon)+ \iota(\stateDimRed; \Upsilon),\label{eqn:err:est:fin}
\end{align}
for $\Upsilon\in\{\switchedSysRedJumps,\switchedSysRedJumpsPert\}$, i.e., either we use the \PBR-\ROM \eqref{eqn:new:red:system} (denoted by $\switchedSysRedJumps$) or the \PBR-\ROM with suitable perturbed Gramians derived in~\Cref{subsec:4.3} (denoted by $\switchedSysRedJumpsPert$).
	\begin{figure}[t]
	\centering
	\begin{subfigure}[t]{.48\linewidth}
		\centering
		\input{Img/Example_2_Plot_1.tex}
		\caption{Qualitative comparison of the first component of the full and reduced output with input signal~$\inp_1$ as in \eqref{eqn:inp:BS} and $20$ switching points.}
	\end{subfigure}
	\hfill
		\begin{subfigure}[t]{.48\linewidth}
		\centering
		\input{Img/Example_2_Plot_2.tex}
		\caption{Qualitative comparison of the second component of the full and reduced outputs with input signal~$\inp_2$ as in \eqref{eqn:inp:BS} and $10$ switching points.}
	\end{subfigure}
	\caption{Example derived from the Black-Scholes model for $\stateDim=1000$. The \GLEs are solved with tolerance $\tol=10^{-8}$.}
		\label{fig4}
\end{figure}

\begin{figure}[t]	
	\centering{
		\begin{subfigure}[t]{0.48\textwidth}
%
\begin{tikzpicture}

\begin{axis}[%
width=\imageWidth,
height=\imageHeight,
at={(1.011in,0.642in)},
scale only axis,
grid=both,
grid style={line width=.1pt, draw=gray!10},
major grid style={line width=.2pt,draw=gray!50},
axis lines*=left,
axis line style={line width=\lineWidth},
mark size=2.2pt,
xmin=0,
xmax=40,
xlabel style={font=\color{white!15!black}},
xlabel={$\stateDimRed$},
ymode=log,
ymin=1e-9,
ymax=1e0,
yminorticks=true,
axis background/.style={fill=white},
legend style={%
	legend cell align=left, 
	align=left, 
	font=\tiny,
	draw=white!15!black,
	at={(0.99,0.99)},
	anchor=north east,},
]
\addplot [color=mycolor1, line width=\lineWidth, mark=o, mark options={solid, mycolor1}]
 table [x=x1, y=y1, col sep=comma]{Img/DataCSV/Example_2_Plot_3.csv};
\addlegendentry{$	\varepsilon(\stateDimRed,\switchedSysRedJumps)$}

\addplot [color=mycolor2, line width=\lineWidth, dashed, mark=o, mark options={solid, mycolor2}]
table [x=x1, y=y2, col sep=comma]{Img/DataCSV/Example_2_Plot_3.csv};
\addlegendentry{	$\varepsilon(\stateDimRed, \switchedSysRedJumpsPert)$}

\addplot [color=mycolor4, line width=\lineWidth, mark=asterisk, mark options={solid, mycolor4}]
table [x=x1, y=y3, col sep=comma]{Img/DataCSV/Example_2_Plot_3.csv};
	\addlegendentry{$\varphi(\stateDimRed, \switchedSysRedJumps)$}
	
\addplot [color=mycolor5, line width=\lineWidth, dashed, mark=asterisk, mark options={solid, mycolor5}]
	table [x=x1, y=y4, col sep=comma]{Img/DataCSV/Example_2_Plot_3.csv};
	\addlegendentry{$ \varphi(\stateDimRed,\switchedSysRedJumpsPert)$}
\end{axis}
\end{tikzpicture}%
			\subcaption{Decay of the error \eqref{eqn:sys:err} and of the error estimate \eqref{eqn:err:est:fin} for the \PBR obtained from $\tilde \fcP(t)$, $\tilde \fcQ(t)$ (see \eqref{eqn:picewise:Gra}) and $\hat \fcP(t)$, $\hat \fcQ(t)$ (see \Cref{subsec:4.3}).}
			\label{fig5:a}
		\end{subfigure}
	\hfill
		\begin{subfigure}[t]{0.48\textwidth}
%
\begin{tikzpicture}

\begin{axis}[%
width=\imageWidth,
height=\imageHeight,
at={(1.011in,0.642in)},
scale only axis,
grid=both,
grid style={line width=.1pt, draw=gray!10},
major grid style={line width=.2pt,draw=gray!50},
axis lines*=left,
axis line style={line width=\lineWidth},
mark size=2.2pt,
xmin=0,
xmax=40,
xlabel style={font=\color{white!15!black}},
xlabel={$\stateDimRed$},
ymode=log,
ymin=5e-9,
ymax=1e0,
yminorticks=true,
axis background/.style={fill=white},
legend style={%
	legend cell align=left, 
	align=left, 
	font=\tiny,
	draw=white!15!black,
	at={(0.99,0.99)},
	anchor=north east,},
]
\addplot [color=mycolor1, line width=\lineWidth, mark=o, mark options={solid, mycolor1}]
 table [x=x1, y=y1, col sep=comma]{Img/DataCSV/Example_2_Plot_4.csv};
\addlegendentry{$\tau(\stateDimRed; \switchedSysRedJumps) $}

\addplot [color=mycolor2, line width=\lineWidth, dashed, mark=o, mark options={solid, mycolor2}]
table [x=x1, y=y2, col sep=comma]{Img/DataCSV/Example_2_Plot_4.csv};
\addlegendentry{$\tau(\stateDimRed; \switchedSysRedJumpsPert)$}

\addplot [color=mycolor4, line width=\lineWidth, mark=asterisk, mark options={solid, mycolor4}]
table [x=x1, y=y3, col sep=comma]{Img/DataCSV/Example_2_Plot_4.csv};
	\addlegendentry{$\chi(\stateDimRed; \switchedSysRedJumps)$}
	
\addplot [color=mycolor3, line width=\lineWidth, dashed, mark=asterisk, mark options={solid, mycolor3}]
	table [x=x1, y=y4, col sep=comma]{Img/DataCSV/Example_2_Plot_4.csv};
	\addlegendentry{$\chi(\stateDimRed; \switchedSysRedJumpsPert)$}
	
\addplot [color=mycolor5, line width=\lineWidth, mark=square, mark options={solid, mycolor5}]
	table [x=x1, y=y5, col sep=comma]{Img/DataCSV/Example_2_Plot_4.csv};
	\addlegendentry{$\iota(\stateDimRed; \switchedSysRedJumps)$}
	
	
\end{axis}
\end{tikzpicture}%
			\subcaption{The components \eqref{eqn:tau:tilde}, \eqref{eqn:chi}, \eqref{eqn:iota} of the error estimate \eqref{eqn:err:est:fin} for the \ROM $\switchedSysRedJumps$ and the \ROM $\switchedSysRedJumpsPert$.}
			\label{fig5:b}
		\end{subfigure}
	}
	\caption{Example derived from the Black-Scholes model for $\stateDim=1000$. \GLEs solved with $\tol=10^{-8}$, input function $\inp_2$ defined in \eqref{eqn:inp:BS}, and $K=10$ switches times.}
	\label{fig5}
\end{figure}
\Cref{fig5:a} shows the actual errors and their corresponding estimators for the two \ROMs. As expected, the error bounds lie above the actual errors in all cases. Notably, the \ROM~ $\switchedSysRedJumpsPert$ obtained using perturbed Gramians exhibits greater accuracy than $\switchedSysRedJumps$, both in terms of the actual error and its estimate when the reduced dimension $\stateDimRed$ is sufficiently large. This behavior is explained in \Cref{fig5:b}, where the error bound is broken down into its three components. For the \PBR-\ROM $\switchedSysRedJumps$, the dominant contribution to the error bound at higher values of $\stateDimRed$ stems from the violation of the \LMI conditions, i.e., the component $\iota(\stateDimRed, \switchedSysRedJumps)$. In contrast, we found that this contribution is negligible (around machine precision) for $\switchedSysRedJumpsPert$, and thus we do not display it. 

\begin{figure}[t]	
	\centering{
		\begin{subfigure}[t]{0.48\textwidth}
%
\begin{tikzpicture}

\begin{axis}[%
width=\imageWidth,
height=\imageHeight,
at={(1.011in,0.642in)},
scale only axis,
grid=both,
grid style={line width=.1pt, draw=gray!10},
major grid style={line width=.2pt,draw=gray!50},
axis lines*=left,
axis line style={line width=\lineWidth},
mark size=2.2pt,
xmin=0,
xmax=40,
xlabel style={font=\color{white!15!black}},
xlabel={$\stateDimRed$},
ymode=log,
ymin=1e-8,
ymax=1e0,
yminorticks=true,
axis background/.style={fill=white},
legend style={%
	legend cell align=left, 
	align=left, 
	font=\tiny,
	draw=white!15!black,
	at={(0.99,0.99)},
	anchor=north east,},
]
\addplot [color=mycolor1, line width=\lineWidth, mark=o, mark options={solid, mycolor1}]
 table [x=x1, y=y1, col sep=comma]{Img/DataCSV/Example_2_Plot_5.csv};
\addlegendentry{$	\varepsilon(\stateDimRed,\switchedSysRedJumpsLE)$}

\addplot [color=mycolor2, line width=\lineWidth, dashed, mark=o, mark options={solid, mycolor2}]
table [x=x1, y=y2, col sep=comma]{Img/DataCSV/Example_2_Plot_5.csv};
\addlegendentry{	$\varepsilon(\stateDimRed, \switchedSysRedJumpsGLE)$}

\addplot [color=mycolor4, line width=\lineWidth, mark=asterisk, mark options={solid, mycolor4}]
table [x=x1, y=y3, col sep=comma]{Img/DataCSV/Example_2_Plot_5.csv};
	\addlegendentry{$\varphi(\stateDimRed, \switchedSysRedJumpsLE)$}
	
\addplot [color=mycolor5, line width=\lineWidth, dashed, mark=asterisk, mark options={solid, mycolor5}]
	table [x=x1, y=y4, col sep=comma]{Img/DataCSV/Example_2_Plot_5.csv};
	\addlegendentry{$ \varphi(\stateDimRed,\switchedSysRedJumpsGLE)$}
\end{axis}
\end{tikzpicture}%
			\subcaption{Decay of the error \eqref{eqn:sys:err} and of the error estimate \eqref{eqn:err:est:fin} for the \PBR obtained via piecewise Lyapunov ($\switchedSysRedJumpsLE$) and piecewise Genaralized Lyapunov ($\switchedSysRedJumpsGLE$) equations.}
			\label{fig6:a}
		\end{subfigure}
	\hfill
		\begin{subfigure}[t]{0.48\textwidth}
%
\begin{tikzpicture}

\begin{axis}[%
width=\imageWidth,
height=\imageHeight,
at={(1.011in,0.642in)},
scale only axis,
grid=both,
grid style={line width=.1pt, draw=gray!10},
major grid style={line width=.2pt,draw=gray!50},
axis lines*=left,
axis line style={line width=\lineWidth},
mark size=2.2pt,
xmin=0,
xmax=40,
xlabel style={font=\color{white!15!black}},
xlabel={$\stateDimRed$},
ymode=log,
ymin=1e-7,
ymax=1e-2,
yminorticks=true,
axis background/.style={fill=white},
legend style={%
	legend cell align=left, 
	align=left, 
	font=\tiny,
	draw=white!15!black,
	at={(0.99,0.99)},
	anchor=north east,},
]


\addplot [color=mycolor4, line width=\lineWidth, mark=asterisk, mark options={solid, mycolor4}]
table [x=x1, y=y3, col sep=comma]{Img/DataCSV/Example_2_Plot_6.csv};
	\addlegendentry{$\chi(\stateDimRed; \switchedSysRedJumpsLE)$}
	
\addplot [color=mycolor3, line width=\lineWidth, dashed, mark=asterisk, mark options={solid, mycolor3}]
	table [x=x1, y=y4, col sep=comma]{Img/DataCSV/Example_2_Plot_6.csv};
	\addlegendentry{$\chi(\stateDimRed; \switchedSysRedJumpsGLE)$}
	
	
	
\end{axis}
\end{tikzpicture}%
			\subcaption{The component \eqref{eqn:chi} of the error estimate \eqref{eqn:err:est:fin} for $\switchedSysRedJumpsLE$ and $\switchedSysRedJumpsGLE$.}
			\label{fig6:b}
		\end{subfigure}
	}
	\caption{Example derived from the Black-Scholes model for $\stateDim=1000$. \GLEs solved with $\tol=10^{-8}$, input function $\inp_2$ defined in \eqref{eqn:inp:BS}, and $K=30$ switches times.}
	\label{fig6}
\end{figure}
We conclude by comparing the \PBR-\ROM obtained from the generalized Lyapunov equations, denoted by $\switchedSysRedJumpsGLE$, with the \PBR-\ROM constructed via the standard Lyapunov equations, denoted by $\switchedSysRedJumpsLE$. As discussed in \Cref{rmk:GLEvsLE}, the latter is expected to yield a larger error due to the increased contribution of the jump terms in the corresponding error bound. To illustrate this effect, we employ the same setup as in \Cref{fig5}, but now restrict attention to a time interval that is $10$ times shorter and introduce $K=30$ switching times, thereby increasing the jump frequency. In \Cref{fig6:a}, we see that the error bound for $\switchedSysRedJumpsGLE$ is several orders of magnitude smaller than that for $\switchedSysRedJumpsLE$; this is also visible in the actual error, which—although not differing as dramatically as the estimate—remains distinctly smaller for $\switchedSysRedJumpsGLE$ when $\stateDimRed$ is large enough. Finally, \Cref{fig6:b} confirms that the discrepancy between the two error estimates is indeed predominantly caused by the jump contribution \eqref{eqn:chi}.

\section{Conclusions and perspectives}

In this manuscript, we study projection-based model reduction for switched linear systems (\SLS) expressed in control form. Specifically, we employ projection matrices obtained through a balancing procedure that is based on approximating the solutions of generalized Lyapunov equations (\GLEs).

We began by introducing a state-of-the-art algorithm for efficiently solving the \GLE. By efficiency, we mean that for large and sparse matrices, the computational cost of computing the approximation scales linearly with respect to $\stateDim$. Our first original contribution is to provide an error bound that allows us to define suitable stopping criteria for a certified approximation of the \GLE solution; that is, given a user-prescribed tolerance $\tol$, the approximation $\tilde \fX$ to the exact solution $\fX$ of the \GLE satisfies $\|\fX - \tilde \fX\|_2 \le \tol$; see \Cref{subsub:gamma:sc}. We describe the conditions under which such a certification can be obtained efficiently; see \Cref{subsubsec:eff:sig:min}. Currently, achieving both efficiency and certification simultaneously for generic cases remains an open problem (even for the simpler Lyapunov equation). One possible direction toward this goal is to rely on the efficient computation of $\varepsilon$-pseudospectra; see the discussion in \Cref{app:A}.

Taking into account the approximation error in the \GLEs solutions is essential for bounding the discrepancy between the output of the original system and that of the reduced model. Under certain restrictions on the system modes, we were able to provide conditions under which the standard error bound obtained through the balancing procedure remains valid; see \Cref{teo4}. Such conditions, however, impose restrictions on the class of linear systems that can be treated; furthermore, the quality of the \ROM approximation may deteriorate if the system modes are ill-conditioned. This motivated us to introduce a novel \MOR procedure that employs projection matrices which are piecewise constant in time, obtained by approximating the solution of as many pairs of \GLEs as there are system modes. We refer to this approach as piecewise-balancing reduction (\PBR) for \SLS; see \Cref{subsec:PBR}. For this procedure, we derived a posteriori error bound (\Cref{teo6}) that takes into account the use of time-varying projection matrices, as well as the possible violation of certain matrix inequalities due to inaccuracies in solving the \GLEs. Computational aspects of \PBR for \SLS and its associated error bound are also discussed in detail in \Cref{subsec:NumAspectsSec4}.

Finally, we carried out numerical experiments to corroborate our theoretical results. Specifically, for a switched system with two dissipative modes, we demonstrate the efficiency and accuracy of \Cref{alg:statIterGLE}, using stopping criteria derived from the error bound established in \Cref{teo1}; see \Cref{fig1}. Subsequently, for a switched system with four modes obtained from $4$ different instances of parameters of a parametric \PDE, we apply the \PBR procedure and evaluate the corresponding error and its error bound; see \Cref{subsec:BlackSholes}.

In future work, we aim to employ the new \PBR methodology in two main directions: first, to construct stable and accurate surrogate models for switched differential-algebraic equations, and second, to speed up model predictive control for \SLS.

\section*{Acknowledgment}
MM acknowledges funding by the BMBF (grant no.~05M22VSA) and acknowledges support by the Stuttgart Center for Simulation Science. BU~is  funded by Deutsche Forschungsgemeinschaft (DFG, German Research Foundation) – Project-ID 258734477 – SFB 1173. Major parts of this manuscript were written while both authors were affiliated with the Stuttgart Center for Simulation Science at the University of Stuttgart.

\bibliographystyle{plain-doi}
\bibliography{journalAbbr,literature}

\appendix
\section{Error bound via \texorpdfstring{$\varepsilon$}{epsilon}-pseudospectrum}\label{app:A}
From \eqref{eqn:errResLyap}, using the integral representation of Lyapunov equations for $\fA$ Hurwitz, one has
      \begin{equation}\label{eqn:starting:bound}
\|\fE_\indKS\|_2\;\le\;\|\fR_\indKS\|_2\int_{0}^{\infty}\|\ee^{t\fA}\|^2_2\,\dt.
      \end{equation}
      To estimate the norm of the matrix exponential, one can resort to the contour integral representation of matrix functions \cite[Cha.~1, Sec.~1.2.3]{Hig08}, i.e.
       \begin{equation}\label{eqn:contour:mat:fun}
\ee^{t\fA}\;=\;\frac{1}{2\pi\mathrm{i}}\int_{\Gamma} \ee^{t z}\left(z\fI_{\stateDim}-\fA\right)^{-1}\,\dz,
      \end{equation}
with $\Gamma$ a closed contour surrounding the singularities of the integrand function, i.e., the eigenvalues of $\fA$. Now, let us consider the $\varepsilon$-pseudospectrum level set, i.e., the set $\Lambda_{\varepsilon}=\{z\,|\,\|(z\fI_{\stateDim}-\fA)^{-1}\|_2\ge\varepsilon^{-1}\}$, and suppose that $\Gamma$ is the circle of radius $R$ and center $z_0$ that contains $\partial \Lambda_{\varepsilon}$ and is such that 
\begin{equation*}
\real(z_0)-R=\alpha_\varepsilon\vcentcolon=\max_{z\in\partial\Lambda_{\varepsilon}}\real(z)<0.   
\end{equation*}
Finally, let us denote by $S_{\varepsilon}$ such a choice of $\Gamma$. Considering the induced Euclidean norm in \eqref{eqn:contour:mat:fun} and using the change of variable $z=z_0+R\ee^{\mathrm{i}\theta}$, we get
\begin{align}
\label{eqn:contour2:0}
\|\ee^{t\fA}\|_2\;\le&\;\frac{1}{2\pi}\int_{S_{\varepsilon}}\ee^{t\real(z)}\|(z\fI_{n}-\fA)^{-1}\|_2\,|\dz|\\
\le&\;\frac{1}{2\pi\varepsilon}\int_{0}^{2\pi}\ee^{t(\real(z_0)+R\cos(\theta))}R\,\dtheta.\label{eqn:contour2}
\end{align}
Substituting \eqref{eqn:contour2} into the r.h.s. of \eqref{eqn:starting:bound}, we find
\begin{align}
\|\fE_{\indKS}\|_2\;\le&\;\frac{R^2\|\fR_{\indKS}\|_2}{4\pi^2\varepsilon^2}\int_{0}^{\infty}\left(\int_{0}^{2\pi}\ee^{t(\real(z_0)+R\cos(\theta))}\,\dtheta\right)^2\dt\label{eqn:bound:1}\\
    \le&\;\frac{R^2\|\fR_{\indKS}\|_2}{2\pi\varepsilon^2}\int_{0}^{\infty}\int_{0}^{2\pi}\ee^{2t(\real(z_0)+R\cos(\theta))}\,\dtheta\,\dt\label{eqn:bound:2}\\
    =&\;\frac{R^2\|\fR_{\indKS}\|_2}{4\pi\varepsilon^2}\int_{0}^{2\pi}-\frac{1}{((\real(z_0)+R\cos(\theta)))}\,\dtheta\label{eqn:bound:3}\\
    =&\;\frac{R^2\|\fR_{\indKS}\|_2}{4\pi\varepsilon^2}\left(\frac{\pi}{\sqrt{\real(z_0)^2-R^2}}\right),\label{eqn:bound:4}
\end{align}
where we apply the Cauchy–Schwarz inequality to pass from \eqref{eqn:bound:1} to \eqref{eqn:bound:2}, use Fubini’s theorem to derive \eqref{eqn:bound:3} from \eqref{eqn:bound:2}, and then perform a standard change of variables in the integral to obtain the final form in \eqref{eqn:bound:4}. Further adjustments give us the bound
\begin{equation}\label{eqn:final:bound}
\|\fE_{\indKS}\|_2\;\le\;C_{\varepsilon}\|\fR_{\indKS}\|_2,\quad\text{with}\quad C_{\varepsilon}\;\vcentcolon=\;\frac{R^2}{4\varepsilon^2\sqrt{|\real(z_0)+R|}\sqrt{|\alpha_{\varepsilon}}|}.
\end{equation}
Therefore, the result of \Cref{prop1} can be replaced by \eqref{eqn:final:bound}. This substitution removes the need to compute $\sigma_{\min}(\fA)$ or other quantities that require $\calO(\stateDim^2)$ floating-point operations. However, it introduces the requirement to compute the $\varepsilon$-pseudospectral abscissa of $\fA$ and to estimate $R$ and $z_0$. These tasks can be computationally expensive and may also yield a large constant $C_{\varepsilon}$, due to the potentially severe overestimation incurred when replacing $\|(z\fI_{\stateDim}-\fA)^{-1}\|_2$ with $\varepsilon^{-1}$. A possible mitigation strategy consists in employing subspace methods for the efficient approximation of the resolvent norm (see \cite{Sir19,ManMG24}), and subsequently performing a numerical quadrature-based approximation of the right-hand side of \eqref{eqn:contour2:0}.
\section{Technical proofs}
\subsection{Proof of \texorpdfstring{\Cref{teo5}}{TEXT}}\label{A-A1}

Consider a switching signal $\switch(t)\in\calS$ and the associated discrete set of times $t_k\in\tilde \calT_{\switch(t)}\vcentcolon=\calT_{\switch(t)}\cup\{t_0,t_{K+1}\}$, where $\calT_{\switch(t)}$ is the set of switching times $t_k$ for $k=1,\ldots,K$ and $t_0\vcentcolon=0$, $t_{K+1}\vcentcolon=t$. Let us define 
\[
\fz(t) \vcentcolon= \fA_{\switch(t_k^{+}),21} \stateRed(t)+ \fB_{\switch(t_k^{+}),2}\inp(t),\quad\quad t\in[t_k,t_{k+1}),\quad k=0,\ldots,K
\] 
where we used the notation introduced in \eqref{eqn:SLS:balanced}, and observe that
\begin{equation}\label{eqn:der:def}
	\dot{\bfX}_o(t)= \bar{\fA}_{\switch(t_k^{+})}\bfX_o(t)+\begin{bmatrix}
		\zeroVec\\
		\fz(t)
	\end{bmatrix},\quad \quad	\dot{\bfX}_c(t)= \bar{\fA}_{\switch(t_k^{+})}\bfX_c(t)-\begin{bmatrix}
		\zeroVec\\
		\fz(t)
	\end{bmatrix}+2\bar{\fB}_{\switch(t_k^{+})}\inp(t),
\end{equation}
with $	\dot{\bfX}_o$ and $\dot{\bfX}_c(t)$ as defined in \eqref{eqn:new:var}. Note that $\bar\fC_{\switch(t)}\bfX_o(t) = \bar \out(t) -\reduce \out(t)$ and for $t\in(t_k,t_{k+1})$ we have
\begin{align}\label{eqn:gen:gra:ine}
	\begin{aligned}
		&\frac{d}{\dt}\left({\bfX}_o^{\T}(t)\fSigma_{\switch(t_k^{+})}{\bfX}_o(t) \right)\;=\;2\bfX^{\T}_o(t)\bar \fA_{\switch(t_k^{+})}^{\T}\fSigma_{\switch(t_k^{+})}\bfX_o(t)+2\sigma({\switch(t_k^{+})})\fz(t)^{\T}\stateBal_2(t)\\
		\le&\;\lambda_{1}(\fM_o(\switch(t_k^{+})))\|\bfX_o(t)\|_2^2-\bfX^{\T}_o(t) \bar \fC^{\T}_{\switch(t_k^{+})}\bar \fC_{\switch(t_k^{+})}\bfX_o(t)+2\sigma({\switch(t_k^{+})})\fz(t)^{\T}\stateBal_2(t),
	\end{aligned}
\end{align}
where we exploited the definition \eqref{eqn:swi:con:fac}, i.e., it holds
\begin{equation*}
\bfX_o(t)^{\T}\left(\bar \fA_{j}^\T\fSigma_{j}+\fSigma_{j}\bar\fA_{j}+\bar \fC^{\T}_{j}\bar \fC_{j}\right)\bfX_o(t)\;\preceq\; \lambda_{1}(\fM_o(j))\|\bfX_o(t)\|_2^2,
\end{equation*}
for all $j\in\switchingSet$. Integrating \eqref{eqn:gen:gra:ine} over the interval $(t_k,t_{k+1})$, we get
\begin{align}\label{eqn:teo4.1:1}
	\begin{aligned}
		\tilde\gleG_{o}&(\switch(t_k^{+}))\;\vcentcolon=\;{\bfX}_o^{\T}(t^{-}_{k+1})\fSigma_{\switch(t^{+}_{k})}{\bfX}_o(t^{-}_{k+1})-{\bfX}_o^{\T}(t^{+}_k)\fSigma_{\switch(t^{+}_k)}{\bfX}_o(t^{+}_k)\\
		=&\;\int_{t_k}^{t_{k+1}} \frac{d}{\dt}\left({\bfX}_o^{\T}(s)\fSigma_{\switch(t_k^{+})}{\bfX}_o(s) \right)\ds\;\le\;\int_{t_k}^{t_{k+1}}\lambda_{1}(\fM_o(\switch(t_k^{+})))\|\bfX_o(s)\|_2^2\,\ds\\
		-&\;\int_{t_k}^{t_{k+1}} \bfX^{\T}_o(s) \bar \fC^{\T}_{\switch(t_k^{+})}\bar \fC_{\switch(t_k^{+})}\bfX_o(s)\,\ds\;+\;\int_{t_k}^{t_{k+1}}2\sigma(\switch(t_k^{+}))\fz(s)^{\T}\stateBal_2(s)\,\ds,
	\end{aligned}
\end{align}
and reordering the terms, we find
\begin{align}\label{eqn:teo4.1:2}
	\begin{aligned}
	\int_{t_k}^{t_{k+1}}\|\bar \out(s)-\tilde \out(s)\|^2_2\,\ds\;\le&\;-\;\tilde \gleG_{o}(\switch(t_k^{+}))+\lambda_{1}(\fM_o(\switch(t_k^{+})))\int_{t_k}^{t_{k+1}}\|\bfX_o(s)\|_2^2\,\ds\;\\
	+&\;2\sigma(\switch(t_k^{+}))	\int_{t_k}^{t_{k+1}}\fz(s)^{\T}\stateBal_2(s)\,\ds.
		\end{aligned}
\end{align}
Next, we use definition \eqref{eqn:new:var} for $\fM_c(j)$, we multiply from left and right by $\fSigma^{-1}_j$ and we find for all $j\in\switchingSet$ that
\begin{equation}\label{eqn:rech:cond}
\fSigma^{-1}_{j}\bar \fA_{j}+\bar \fA_{j}^\T\fSigma^{-1}_{j}+\fSigma^{-1}_{j}\bar \fB_{j}\bar \fB_{j}^{\T}\fSigma^{-1}_{j}-\lambda_{1}(\fM_c(j))\fSigma^{-2}_{j}\;\preceq \;\zeroMat.
\end{equation}
Using the Schur complement it follows that \eqref{eqn:rech:cond} is equivalent to
\begin{align*}
\begin{aligned}
\begin{bmatrix}
\fSigma^{-1}_{j}\bar \fA_{j}+\bar \fA_{j}^\T\fSigma^{-1}_{j}-\lambda_{1}(\fM_c(j))\fSigma^{-2}_{j}&\fSigma^{-1}_{j}\bar \fB_{j}\\
\bar \fB_{j}^{\T}\fSigma^{-1}_{j}&-\fI
\end{bmatrix}\preceq\zeroMat,
\end{aligned}
\end{align*}
then,  by definition \eqref{eqn:der:def}, it follows for $t\in(t_k,t_{k+1})$ that
\begin{align*}
\begin{aligned}
&\frac{d}{\dt}\left({\bfX}_c^{\T}(t)\fSigma^{-1}_{\switch(t_k^{+})}{\bfX}_c(t) \right)\;\\
=&\;\begin{bmatrix}
{\bfX}_c(t)\\
2\inp(t)
\end{bmatrix}^{\T}\begin{bmatrix}
\fSigma^{-1}_{\switch(t_k^{+}))}\bar \fA_{\switch(t_k^{+}))}+\bar \fA_{\switch(t_k^{+}))}^\T\fSigma^{-1}_{\switch(t_k^{+}))}-\lambda_{1}(\fM_c(\switch(t_k^{+}))))\fSigma^{-2}_{\switch(t_k^{+}))}&\fSigma^{-1}_{\switch(t_k^{+}))}\bar \fB_{\switch(t_k^{+}))}\\
\bar \fB_{\switch(t_k^{+}))}^{\T}\fSigma^{-1}_{\switch(t_k^{+}))}&-\fI
\end{bmatrix}\begin{bmatrix}
{\bfX}_c(t)\\
2\inp(t)
\end{bmatrix}\\
+&\;4\|\inp(t)\|_2^2
\;+\;\lambda_{1}(\fM_c(\switch(t_k^{+})))\bfX_c(t)^{\T}\fSigma^{-2}_{\switch(t_k^{+})}\bfX_c(t)\;-\;2\bfX_c(t)\fSigma^{-1}_{\switch(t_k^{+})}\begin{bmatrix}\zeroVec\\
	\fz(t)
\end{bmatrix}\\
\le &\;\lambda_{1}(\fM_c(\switch(t_k^{+})))\frac{1}{\sigma(\switch(t_k^{+}))^{2}}\|\bfX_c(t)\|_2^2\;+\;4\|\inp(t)\|_2^2-2\frac{1}{\sigma(\switch(t_k^{+}))}\fz(t)^{\T}\stateBal_2(t)
\end{aligned}
\end{align*}
where we also used the fact that $\bfX_c(t)\fSigma^{-1}_{\switch(t_k^{+})}\begin{bmatrix}\zeroVec&
	\fz(t)
\end{bmatrix}^{\T}={\sigma(\switch(t_k^{+}))^{-1}}\fz(t)^{\T}\stateBal_2(t)$. Now we integrate over the interval $(t_k,t_{k+1})$ and we get
\begin{align*}
	\begin{aligned}
		\tilde\gleG_{c}(\switch(t_k^{+}))\;\vcentcolon=&\;{\bfX}_c^{\T}(t^{-}_{k+1})\fSigma^{-1}_{\switch(t^{+}_{k})}{\bfX}_c(t^{-}_k)-{\bfX}_c^{\T}(t^{+}_k)\fSigma^{-1}_{\switch(t^{+}_k)}{\bfX}_c(t^{+}_k)   
		=\int_{t_k}^{t_{k+1}} \frac{d}{\dt}\left({\bfX}_c^{\T}(s)\fSigma^{-1}_{\switch(t_k^{+})}{\bfX}_c(s) \right)\ds\\
		\le&\;4\int_{t_k}^{t_{k+1}}\|\inp(s)\|_2^2\,\ds-2\int_{t_k}^{t_{k+1}}\frac{1}{\sigma(\switch(t_k^{+}))}\fz(s)^{\T}\stateBal_2(s)\,\ds\\
		+&\;\lambda_{1}(\fM_c(\switch(t_k^{+})))\frac{1}{\sigma(\switch(t_k^{+}))^{2}}\int_{t_k}^{t_{k+1}}\|\bfX_c(s)\|_2^2\,\ds,
	\end{aligned}
\end{align*}
and after rearranging the terms, we have the relation
\begin{align}\label{eqn:rec:imp:rel}
	\begin{aligned}
		2\int_{t_k}^{t_{k+1}}\fz(s)^{\T}\stateBal_2(s)\,\ds\;\le&\;-\tilde \sigma(\switch(t_k^{+}))\tilde \gleG_{c}(\switch(t_k^{+}))+4 \sigma(\switch(t_k^{+})) \int_{t_k}^{t_{k+1}}\|\inp(s)\|_2^2\,\ds\\
		+&\;\lambda_{1}(\fM_c(\switch(t_k^{+})))\frac{1}{\sigma(\switch(t_k^{+}))}\int_{t_k}^{t_{k+1}}\|\bfX_c(s)\|_2^2\,\ds.
	\end{aligned}
\end{align}
Multiplying \eqref{eqn:rec:imp:rel} by $\sigma(\switch(t_k^{+}))$ and then substituting it into \eqref{eqn:teo4.1:2} gives us
\begin{align}\label{eqn:fun:rel:rub:int}
	\begin{aligned}
		&\int_{t_k}^{t_{k+1}}\|\bar \out(s)-\tilde \out(s)\|^2_2\,\ds
		\;\le\;4\sigma(\switch(t_k^{+}))^{2}\int_{t_k}^{t_{k+1}}\|\inp(s)\|_2^2\,\ds\\
		-&\;\tilde \gleG_{o}(\switch(t_k^{+}))-\sigma(\switch(t_k^{+}))^{2}\tilde \gleG_{c}(\switch(t_k^{+}))\\
		+&\;\lambda_{1}(\fM_o(\switch(t_k^{+})))\int_{t_k}^{t_{k+1}}\|\bfX_o(s)\|_2^2\,\ds+\lambda_{1}(\fM_c(\switch(t_k^{+})))\int_{t_k}^{t_{k+1}}\|\bfX_c(s)\|_2^2\,\ds.
	\end{aligned}
\end{align}
Finally, repeating the previous passages for all $k=0,\ldots,K$, summing up all the relations of type \eqref{eqn:fun:rel:rub:int}, and recalling the definition \eqref{eqn:swi:con:fac} for $\gleG_{o,\switch}(\stateDim,\stateDimRed)$, $\gleG_{c,\switch}(\stateDim,\stateDimRed)$, $\gleH_{o,\switch}(\stateDim,\stateDimRed)$, and $\gleH_{c,\switch}(\stateDim,\stateDimRed)$ we end with
\begin{align*}
	\begin{aligned}
		\int_{0}^{t}\|\bar \out(s)-\tilde \out(s)\|^2_2\,\ds
		\le
		\gleH_{o,\switch}(\stateDim,\stateDimRed)+\gleH_{c,\switch}(\stateDim,\stateDimRed)-\gleG_{o,\switch}(\stateDim,\stateDimRed)-\gleG_{c,\switch}(\stateDim,\stateDimRed)
		+4\tilde \sigma^{2}\int_{0}^{t}\|\inp(s)\|_2^2\,\ds,
	\end{aligned}
\end{align*}
that corresponds to statement \eqref{eq:err:est:2} by noting that $\out\equiv\bar \out$.
\end{document}